\documentstyle[12pt,fleqn,leqno]{article}
\newtheorem{thm}{Theorem}[section]
\newcommand{\be}{\begin{equation}}
\newcommand{\ee}{\end{equation}}
\newcommand{\ba}{\begin{array}}
\newcommand{\ea}{\end{array}}

\newcommand{\dq}{\delta_{q}}
\renewcommand{\a}{\alpha}
\renewcommand{\b}{\beta}
\renewcommand{\l}{\lambda}
\renewcommand{\t}{\theta}

\newcommand{\mg}{\rm}
\renewcommand{\em}{\it}
\newcommand{\bea}{\begin{eqnarray}}
\newcommand{\eea}{\end{eqnarray}}

\newcommand{\dis}{\vspace{-0.6\abovedisplayskip}}
\newcommand{\Int}{\int_{-1}^1}
\newcommand{\Sum}{\sum_{n=1}^\infty}
\newcommand{\Summ}{\sum_{n=0}^\infty}
\newcommand{\Prod}{\prod_{n=0}^\infty}

\newcommand{\g}{\gamma}
\newcommand{\G}{\Gamma}
\begin{document}
\newtheorem{lem}[thm]{Lemma}
\newtheorem{cor}[thm]{Corollary}
\title{\bf Diagonalization of Certain Integral Operators II\thanks{
Research partially supported  by NSF grant DMS 9203659 and NSERC
grant 
A6197.} }
\author{Mourad E. H.  Ismail, Mizan Rahman, and Ruiming Zhang}
\date{}
\maketitle
\begin{abstract}
We establish an integral representations of a right inverses of the
Askey-Wilson 
finite difference operator in an $L^2$ space  weighted by the
weight function of
the continuous $q$-Jacobi polynomials.  We  characterize  the
eigenvalues 
of this integral operator and prove a $q$-analog of the expansion
of $e^{ixy}$
in Jacobi polynomials of argument $x$. We also outline a general
procedure
of finding integral representations for inverses of linear
operators.
\end{abstract}

\bigskip 
{\bf Running title}:$\;$ Integral Operators


\vfill\eject
\section{Introduction.} 
The Askey-Wilson divided difference operator ${\cal D}_q$,
\cite{As:Wi} is defined by
\be
({\cal D}_q f)(x)=\frac{\dq\,f(x)}{i(q^{1/2}-q^{-1/2})\sin\t},\quad
x=\cos\t,
\ee
where 
\be
(\dq g)(e^{i\theta})=g(q^{1/2}e^{i\t })-g(q^{-1/2}e^{i\t}).
\label{a2}
\ee
Observe that $i(q^{1/2}-q^{-1/2})\sin\t$ which appears in the
denominator
of (1.1) is $\dq x$, $x$ being the identity map: $x\mapsto x$
evaluated at
$x$.  Magnus \cite{Ma} showed how the Askey-Wilson operator arises 
naturally from divided difference operators.

It is easy to that 
\bea
{\cal D}_qT_n(x) = \frac{q^{n/2} - q^{-n/2}}{q^{1/2} -
q^{-1/2}}U_{n-1}(x), \nonumber
\eea
where $T_n(x)$ and $U_n(x)$ are Chebyshev polynomials of the 
first and second kinds, respectively. 
Therefore ${\cal D}_q$ maps a polynomial of degree $n$ to a
polynomial of degree $n-1$. As such ${\cal D}_q$ resembles the 
differential operator.
It was observed in \cite{Br:Is} and \cite{Is:Zh} that one can
construct integral operators
which are a right inverse to ${\cal D}_q$ on certain weighted
spaces. In
\cite{Is:Zh} Ismail and Zhang diagonalized the right inverse to
${\cal D}_q$ on
weighted $L_2$ space on $[-1,1]$ with Jacobi weights
$(1-x)^\a(1+x)^\b$ or
$q$-ultraspherical weights
\[
w_{\b;q}(x)=\frac1{\sqrt{1-x^2}}\Prod\frac{1-2(2x^2-1)q^n+q^{2n}}
{1-2(2x^2-1)q^{n+\nu}+q^{2n+2\nu}},\quad \b=q^\nu,\;\nu>0.
\]

Recall the notations
\be
(a;\;q)_0:=1,\quad (a;\;q)_n:=\prod_{j=1}^{n}(1-aq^{j}),\quad
n=1,2,\ldots,
\mbox{ or }\infty,
\ee
\be
(a_1,\ldots,a_m;\:q)_n=\prod_{k=1}^{m}\,(a_k;\;q)_n,
\ee
for $q$-shifted factorials.
We shall normally drop ``$;q$'' from the shifted factorials in
(1.3) and
(1.4) when this does not lead to any confusion. Thus for the
purpose of
this paper we will use
\be
(a)_n:=(a;q)_n
\ee
A basic hypergeometric series is
\bea\lefteqn{
{}_r\phi_s\left(\left.\ba{c} a_1,\ldots,a_r \\ b_1,\ldots,b_s \ea
\right|\,q,\;z
\right)=\;_r\phi_s(a_1,\dots,a_r;b_1,\dots,b_s;q,z)}\\
& & :=\sum_{n=0}^\infty
\frac{(a_1,\ldots,a_r;q)_n}{(q,b_1,\ldots,b_s;q)_n}\,z^n\,
\left[(-1)^nq^{n(n-1)/2}\right]^{1+s-r}.\nonumber
\eea

A $q$-analog of Jacobi polynomials is given by
\be
P_n^{(\a,\b)}(x;q)=\frac{(q^{\a+1},-q^{\b+1};q)_n}{(q,-q;q)_n}
\mbox{}_4\phi_3\left(\left.\ba{c}
q^{-n},q^{n+\a+\b+1},q^{1/2}e^{i\t},q^{1/2}e^{-i\t}\\
\quad q^{\a+1},\quad -q^{\b+1},\quad -q \ea \right|q,q\right),
\ee
where $x=\cos \t$. The $P_n$'s are called the continuous $q$-Jacobi
polynomials. 
In what follows we assume
\be
e^{i\t}=x+\sqrt{x^2-1},\quad
e^{-i\t}=x-\sqrt{x^2-1},
\ee
where the sign of the square root in (1.8) is taken so that
$\sqrt{x^2-1}
\approx x$ as $x\to \infty$ in the complex plane. The normalization
in
(1.7) was introduced by Rahman in \cite{Ra}. The original
normalization  
used by Askey and Wilson in \cite{As:Wi} is $P_n^{(\a,\b)}(x|q)$,
where
\be
P_n^{(\a,\b)}(x;\:q)=\frac{(-q^{\a+\b+1};q)_n}{(-q;q)_n}q^{-\a
n}P_n^{(\a,\b)}(x|q^2).
\ee
Both $\{P_n^{(\a,\b)}(x;q)\}$ and $\{P_n^{(\a,\b)}(x|q)\}$ are
called 
continuous $q$-Jacobi polynomials because they are orthogonal with
respect
to an absolutely continuous measure. The orthogonality relation is
\be
\int_{-1}^1P_n^{(\a,\b)}(x|q)P_m^{(\a,\b)}(x|q)w_{\a,\b}(x|q)\,dx=
h_n^{(\a,\b)}(q)\delta_{m,n}.
\ee
The weight function $w_{\a,\b}(x|q)$ and the normalization
constants
$h_n^{(\a,\b)}(q)$ are given by \cite[(7.5.28), (7.5.30),
(7.5.31)]{Ga:Ra}.
\be w_{\a,\b}(\cos\t|q^2)=
\frac{(e^{2i\t},e^{-2i\t};q^2)_\infty\,(1-x^2)^{-1/2}}{(q^{\a+1/2
}e^{i\t},
q^{\a+1/2}e^{-i\t}, -q^{\b+1/2}e^{i\t},
-q^{\b+1/2}e^{-i\t};q)_\infty},
\ee
\dis\dis
\be
h_n^{(\a,\b)}(q)=\frac{2\pi(1-q^{\a+\b+1})(q^{(\a+\b+2)/2},q^{(\a
+\b+3)/2};
q)_\infty}{(q,q^{\a+1},q^{\b+1},-q^{(\a+\b+1)/2},-q^{(\a+\b+2)/2}
;q)_\infty}
\ee
\dis\dis
\[\mbox{\hspace{1.1in}}
\cdot\frac{(q^{\a+1},q^{\b+1},-q^{(\a+\b+3)/2};q)_n\,
q^{n(2\a+1)/2}}{
(1-q^{2n+\a+\b+1})(q,q^{\a+\b+1},-q^{(\a+\b+1)/2};q)_n}.
\]
Note that (7.5.31) in \cite{Ga:Ra} contains a misprint where
"$q^{n(2\a+1)/4}$"
on the right-hand side should read "$q^{n(2\a+1)/2}$".

The Askey-Wilson operator acts on continuous $q$-Jacobi polynomials
in a
very natural way. It's action is given by, \cite[ (7.7.7)]{Ga:Ra},
\be
{\cal D}_q
P_n^{(\a,\b)}(x|q)=\frac{2q^{-n+(2\a+5)/4}(1-q^{\a+\b+n+1})}{
(1+q^{(\a+\b+1)/2})(1+q^{(\a+\b+2)/2})(1-q)}P_{n-1}^{(\a+1,\b+1)}
(x|q).
\ee
Following \cite{Is:Zh} we define ${\cal D}_q$ densely on
$L_2[w_{\a,\b}(x|q)]$ by
\be
{\cal D}_q f \sim \Sum \frac{2q^{-n+(2\a+5)/4}(1-q^{\a+\b+n+1})}{
(1+q^{(\a+\b+1)/2})(1+q^{(\a+\b+2)/2})(1-q)}f_n
P_{n-1}^{(\a+1,\b+1)}(x|q),
\ee
\dis
if
\dis
\be
f\sim \Summ f_n P_{n}^{(\a,\b)}(x|q).
\ee
Clearly ${\cal D}_q$ as defined by (1.14) and (1.15) maps a dense
subset 
of $L_2[w_{\a,\b}(x|q)]$ into  $L_2[w_{\a+1,\b+1}(x|q)]$.
We are interested in finding a formal inverse to ${\cal D}_q$, that
is we
seek a linear operator $T_{\a,\b;q}$
 which maps  $L_2[w_{\a+1,\b+1}(x|q)]$ into $L_2[w_{\a,\b}(x|q)]$
such
that ${\cal D}_q\;T_{\a,\b;q}$ is the identity map on the range of
${\cal D}_q$. It is clear from (1.14) and (1.15) that we may
require 
$T_{\a,\b;q}$ to satisfy
\be
(T_{\a,\b;q}\:g)(x)\sim \Summ
\frac{(1-q)(-q^{(\a+\b+1)/2};q^{1/2})_2
q^{n-(2\a+1)/4}}{2(1-q^{\a+\b+n+2})}g_n P_{n+1}^{(\a,\b)}(x|q)
\ee
for
\be
g(x)\sim \Summ g_n P_{n}^{(\a+1,\b+1)}(x|q).
\ee
It is easy to find a representation of $T_{\a,\b;q}$ as an integral

operator. We use the orthogonality relation (1.10) to write $g_n$
as
\[
\Int g(x)w_{\a+1,\b+1}(x|q)P_{n}^{(\a+1,\b+1)}(x|q)\;
dx/h_n^{(\a+1,\b+1)}(q)
\]
and formally interchange summation and integration in (1.16). The
result
is the formal definition
\be
(T_{\a,\b;q}\:g)(x)=\Int K_{\a,\b;q}(x,y)\; g(y) \;
w_{\a+1,\b+1}(y|q)dy,
\ee
and the kernel $K_{\a,\b;q}(x,y)$ is given by
\be
 K_{\a,\b;q}(x,y) =\Summ
\frac{(1-q)(-q^{(\a+\b+1)/2};q^{1/2})_2
}{2(1-q^{\a+\b+n+2}) h_n^{(\a+1,\b+1)}(q)}\; q^{n-(2\a+1)/4}
P_{n+1}^{(\a,\b)}(x|q)
P_{n}^{(\a+1,\b+1)}(y|q)
\ee

The purpose of this work is to study the spectral properties of the
integral operator (1.18). It is worth noting that $T_{\a,\b;q}$ is
linear
but is not normal and a spectral theory of such operators is not
readily
available. Our main result is Theorem 3.1 which characterizes the
eigenvalues
of $T_{\a,\b;q}$ as zeros of a certain transcendental
function. In order to prove Theorem 3.1 we proved several auxiliary
results
which may be of interest by themselves. First in \S2 we solve the
connection coefficient problem of expressing $P_n^{(\a,\b)}(x)$ in
terms of
$\{P_j^{(\a+1,\b+1)}(x)\}$.
The solution of this connection coefficient problem is then used to
expand
$w_{\a+1,\b+1}(x|q)P_{n}^{(\a+1,\b+1)}(x|q)$ in terms of
$\{w_{\a,\b}(x|q)
P_{k}^{(\a,\b)}(x)\}$. In Section 3 we used the latter expansion to
find a
tridiagonal matrix representation of $T_{\a,\b;q}$. In Section 3 we
also
find the eigenvalues and eigenfunctions $g(x|\l)$ so that
\[
T_{\a,\b;q}\:g =\l g.
\]
The eigenvalues are multiples of the reciprocals of the zeros of a
transcendental function  $X_{-1}^{(\a,\b)}(1/x)$ defined in (5.13).
Such a function is a $q$-analog of the confluent hypergeometric
function $_1F_1$. The
 eigenfunctions are shown to have the orthogonal expansion
\be
\sum_1^{\infty} a_n(\lambda|q) P_n^{(\a, \b)}(x;q)
\ee
such that
\be
\sum_1^{\infty} h_n^{(\a, \b)}(q)\,|a_n(\lambda|q)|^2\, < \,
\infty.
\ee
We may normalize the eigenfunctions by choosing $a_1(\lambda|q) =
1$. 
With this normalization we prove , in Section 3, that
$a_n(\lambda|q)$ 
is a polynomial of degree $n$. The $a_n$'s are q-analogs of
polynomials
studied by  Walter Gautschi \cite{Gau}  and Jet Wimp \cite{Wi}. An
explicit formula for $a_n(\lambda|q)$ 
is given in Section 4, see (3.11) and (4.1). The large $n$
asymptotics of 
$a_n(\lambda|q)$  in different parts of the complex plane are found
in
 Section 5 and they are used to characterize the $\lambda$'s for
which (1.21) holds.

In Section 5 we note that when $\a$ and $\b$ are complex conjugates
and are not real  then the polynomials 
$i^{-n} a_n(ix|q)$ are real orthogonal polynomials.  
They are orthogonal with respect to a discrete measure supported at
the zeros of a transcendental function which we denoted by
$F_L(\eta,\rho;q)$. The function $F_L(\eta,\rho;q)$
is a $q$-analog of the regular Coulomb wave function
$F_L(\eta,\rho)$,
\cite[Chapter 14]{Ab:St}. Our analysis implies that the functions
$F_L(\eta,\rho;q)$
have only real and simple zeros. Neither the functions
$F_L(\eta,\rho;q)$
nor their zeros seem to have been studied before this work.

In Section 6 we prove that the eigenfunctions are constant
multiplies of the $q$-exponential function
${\cal E}_q(x;-i,b)$, where
\be 
{\cal E}_q(x;a,b):=\Summ
\frac{q^{n^2/4}}{(q)_n}(aq^{(1-n)/2}e^{i\t},
 aq^{(1-n)/2}e^{-i\t})_n b^n.
\ee
The function ${\cal E}_q$ was introduced in \cite{Is:Zh}. Since 
the eigenfunctions also have the orthogonal expansion (1.20) we
obtain
 an identity valid on the spectrum of $T_{\a, \b;q}$. This identity
is (6.13)
and is a $q$-analog of 
\be
e^{ixy}= e^{-iy} \sum_0^{\infty} \, \frac{(\a + \b +1)_n}{(\a + \b
+1)_{2n}}
 \; (2iy)^n{}_1F_1\left(\left.\ba{c} n +\b+1 \\
2n +\a + \b + 2 \ea \right|2iy\right) \; P_n^{(\a,\b)}(x),\quad -
1 < x < 1,
\ee
see (10.20.4) in \cite{Er:Ma}.
In Section 6 we use properties of basic hypergeometric functions to
show 
 that the above mentioned $q$-identity holds also off the spectrum 
of  $T_{\a, \b;q}$.  In Section 7 we give a second proof  of (6.13)
using a technique
similar to what was used in \cite{Is:Zh} to prove the same result
for the 
continuous $q$-ultraspherical polynomials. We also include in \S 7 
a 
formal approach  to finding the spectrum of certain integral
operators of the
type considered in this paper. 
In \S 8 we include  some remarks on  asymptotic results of Schwartz
\cite{Sc}, Dickinson, Pollack and Wannier 
\cite{Di:Po} and general remarks on this work.

In many of our calculations we found it advantageous to follow
\cite{Ga:Ra}.
The only disadvantage is that we have to introduce some additional
relations. Recall that the Askey-Wilson polynomials are defined by
\cite[(7.5.2)]{Ga:Ra}
\be
p_n(x;a,b,c,d |q) =(ab,ac,ad;q)_n a^{-n}
{}_4\phi_3\left(\left.\ba{c} q^{-n}, abcd q^{n-1}, ae^{i\t},
ae^{-i\t} \\
ab,\quad ac,\quad ad \ea \right|q;q\right).
\ee
Their orthogonality relation is \cite[(7.5.15), (7.5.16)]{Ga:Ra},
\cite{As:Wi}
\be
\int_{-1}^{1} \frac{h(x;1, -1, q^{1/2}, -q^{1/2})}{h(x; a, b, c,
d)}
p_n(x; a, b, c, d)p_m(x; a, b, c, d) \frac{dx}{\sqrt{1-x^2}}
\ee
\bea
= \kappa(a,b,c,d|q) \frac{(1-abcd/q)(q, ab, ac, ad, bc, bd, cd)_n}
{(1 - abcdq^{2n-1})(abcd/q)_n} \delta_{m,n}, \nonumber
\eea
with
\be
h(\cos \theta; a_1, a_2, a_3, a_4):= \prod_{j=1}^{4}(a_j
e^{i\theta}, a_j e^{-i\theta})_{\infty}, 
\ee
and
\be
\kappa(a,b,c,d|q)= 2\pi
(abcd)_\infty/[(q,ab,ac,ad,bc,bd,cd)_\infty].
\ee
We shall also use the notation
\be 
h(\cos \theta; a) := (a e^{i\theta}, a e^{-i\theta})_{\infty}.
\ee
Note that
\be 
h(\cos \theta; 1, -1, q^{1/2}, q^{-1/2}) = (e^{2i\theta}, 
e^{-2i\theta})_{\infty}.
\ee
Note also that the continuous $q$-Jacobi polynomials correspond to
the identification of parameters 
\be
a = q^{(2\a + 1)/4}, \quad b= q^{(2\a + 3)/4}, \quad c = -q^{(2\b
+ 1)/4},
 \quad  d  =- q^{(2\b + 3)/4}.
\ee
In fact the polynomials $P_n^{(\a, \b)}(x|q)$  of (1.7) 
have the alternate representation, 
\be
P_n^{(\a, \b)}(x|q) = \frac{(q^{\a+1};q)_n}{(q; q)_n} \;
{}_4\phi_3\left(\left.\ba{c} q^{-n},q^{n +\a + \b + 1},
q^{(2\a+1)/4} e^{i\t}, q^{(2\a+1)/4} e^{-i\t} \\
q^{\a + 1},\; \;  -q^{(\a + \b + 1)/2},\; \;  -q^{(\a + \b + 2)/2} 
\ea \right|q;q\right).
\ee

\newpage
\section{Connection Coefficients}
\setcounter{equation}{0}
In this section we derive  a $q$-analogue of the formula
\be
(1-x^2)P_{n-1}^{(\a+1,\b+1)}(x)=\frac{4(n+\a)(n+\b)}{(2n+\a+\b)(2
n+\a+\b+1)}
P_{n-1}^{(\a,\b)}(x)
\ee
\dis\dis
\[
+\frac{4n(\a -\b)}{(2n+\a+\b)(2n+\a+\b+2)}P_{n}^{(\a,\b)}(x)-
\frac{4n(n+1)}{(2n+\a+\b+1)(2n+\a+\b+2)}P_{n+1}^{(\a,\b)}(x).
\]
Formula (2.1) is (1) in \S37 of \cite{Rai}. The Jacobi polynomials
$\{P_n^{(\a,\b)}(x)
\}$ are orthogonal with respect to  $(1-x)^\a(1+x)^\b$ and (2.1) is
essentially the expansion of  $(1-x)^{\a+1}(1+x)^{\b+1}
P_{n-1}^{(\a+1,\b+1)}(x)$ in terms of
$\{(1-x)^\a(1+x)^\b P_j^{(\a,\b)}(x)\}_{j=0}^\infty$. The
$q$-analog of this
question is to expand $w_{\a+1,\b+1}(x;q)
P_{n-1}^{(\a+1,\b+1)}(x;q)$ in terms of
$\{w_{\a,\b}(x;q) P_{j}^{(\a,\b)}(x;q)\}_{j=0}^\infty$. 
It is worth mentioning that the latter problem is equivalent to the
expression of $P_n^{(\a,\b)}(x|q)$ in terms of
$\{P_j^{(\a+1,\b+1)}(x|q)\}_{j=0}^n$.
This follows from the following known observation. If $\{
p_n(x;\l)\}$ are
orthonormal with respect to a weight function $w(x;\l)$ then the 
connection coefficient formula
\[
p_n(x;\l)=\sum_{j=0}^n c_{n,j}(\l,\mu) p_j(x;\mu)
\]
 holds if and only if its dual, namely
\[
w(x;\mu) p_n(x;\mu)\sim \sum_{j=n}^\infty c_{j,n}(\l,\mu)
w(x;\l) p_j(x;\l),
\]
holds.  This latter fact follows from computing the Fourier
coefficients of both
sides.

The main result of this
section is the following theorem

\begin{thm}
We have
\bea\lefteqn{
(1-2xq^{\a+1/2}+q^{2\a+1})(1+2xq^{\b+1/2}+q^{2\b+1})P_{n-1}^{(\a+
1,\b+1)}(x;q)
}\\
& & =\frac{(1+q^{\a+\b+n})(1+q^{\a+\b+n+1})(1+q^{\a+n})(1+q^{\b+n})
(1-q^{\a+n})(1-q^{\b+n})}{(1-q^{2n+\a+\b})(1-q^{2n+\a+\b+1})}P_{n
-1}^{(\a,\b)}(x;q)
\nonumber \\
& & +\frac{(1+q^{\a+\b+n+1})(1+q^{\a+\b+2n+1})(1+q^{n})^2
(1-q^{n})(1-q^{\a-\b})}{(1-q^{2n+\a+\b})(1-q^{2n+\a+\b+2})}q^{\b}
P_{n}^{(\a,\b)}(x;q) \nonumber\\
& & -\frac{(1+q^{n})^2(1+q^{n+1})^2(1-q^{n})
(1-q^{n+1})}{(1-q^{2n+\a+\b+1})(1-q^{2n+\a+\b+2})}q^{\a+\b}
P_{n+1}^{(\a,\b)}(x;q). \nonumber
\eea
\end{thm}

For our purposes it is convenient to express (2.2) in the
Askey-Wilson
normalization $P_n^{(\a,\b)}(x|q)$. The result is 
\bea\lefteqn{
(1-2xq^{\a+1/2}+q^{2\a+1})(1+2xq^{\b+1/2}+q^{2\b+1})P_{n-1}^{(\a+
1,\b+1)}(x|q^2)
}\\
& & =\frac{(1-q^{2\a+2n})(1-q^{2\b+2n})(-q^{\a+\b+1};q)_2}{
(q^{2n+\a+\b};q)_2}q^{n-1}P_{n-1}^{(\a,\b)}(x|q^2) \nonumber \\
& &
+\frac{(-q^{\a+\b+1};q)_2(1+q^{\a+\b+2n+1})(1-q^{2n})(1-q^{\a-\b}
)}{
(q^{2n+\a+\b};q^2)_2}q^{\b-\a+n-1}P_{n}^{(\a,\b)}(x|q^2) \nonumber
\\
& & -\frac{(-q^{\a+\b+1};q)_2(1-q^{2n})(1-q^{2n+2})}{
(q^{2n+\a+\b+1};q)_2}q^{\b-\a+n-1}P_{n+1}^{(\a,\b)}(x|q^2).
\nonumber
\eea

The rest of this section will be devoted to proving (2.2). Our
proof is
very technical and the reader who is more conceptually oriented is
advised
to turn to Section 3.

Our proof uses the Sears transformation \cite[(III.15)]{Ga:Ra}
\be
{}_4\phi_3\left(\left. \ba{c} q^{-n}, a, b, c\\ d, e, f
\ea\right|q,q\right)=
\frac{(e/a, f/a)_n}{(e,f)_n}\,a^n\,
{}_4\phi_3\left(\left. \ba{c} q^{-n}, a, d/b, d/c\\ d, aq^{1-n}/e,
aq^{1-n}/f 
\ea\right|q,q\right),
\ee
where $def=abcq^{1-n}$. Note that the parameters $q^{-n}, a$ and
$d$  remain
invariant under the transformation (2.4).

We now proceed with the proof. We seek a connection coefficient
formula of
the type
\bea\lefteqn{
(1-2b\cos \t+b^2)(1+2c\cos \t +c^2){}_4\phi_3\left(\left.\ba{c}
q^{1-n},bcq^{n+1},q^{1/2}e^{i\t},q^{1/2}e^{-i\t}\\
bq^{3/2},-cq^{3/2},-q \ea\right|q,q\right)} \nonumber\\
& &=\sum_{k=0}^{n+1} A_k\:{}_4\phi_3\left(\left.\ba{c}
q^{-k},bcq^{k},q^{1/2}e^{i\t},q^{1/2}e^{-i\t}\nonumber\\
bq^{1/2},-cq^{1/2},-q \ea\right|q,q\right), 
\eea
where $b:=q^{\a+1/2},\; c:=q^{\b+1/2}$.

The orthogonality relation (1.23) gives
\bea\lefteqn{
A_k\,\kappa(q^{1/2}, b, -c,-q^{1/2})\frac{(1-bc)(q,-bc,-b\sqrt{q},
c\sqrt{q})_k}{(1-bcq^{2k})(bc,-q,-c\sqrt{q}, b\sqrt{q})_k} q^k
\nonumber}\\
& & =\int_0^\pi\frac{(e^{2i\t},e^{-2i\t})_\infty}{
h(\cos\t;\sqrt{q}, bq, -cq, -\sqrt{q})} \nonumber \\
& & \mbox{\hspace{0.3in}} \cdot{}_4\phi_3 \left(\left. \ba{c}
q^{-k}, bcq^k,
q^{1/2}e^{i\t},q^{1/2}e^{-i\t} \\ bq^{1/2},\qquad -cq^{1/2},\qquad
-q \ea
\right|q,q\right) \nonumber \\
& & \mbox{\hspace{0.3in}} \cdot{}_4\phi_3 \left(\left. \ba{c}
q^{1-n}, bcq^{n+1},
q^{1/2}e^{i\t},q^{1/2}e^{-i\t} \\ bq^{3/2},\qquad -cq^{3/2},\qquad
-q \ea
\right|q,q\right)\, d\t. \nonumber 
\eea
Apply the Sears transformation (2.4) with invariant parameters
$q^{-k}, bcq^k$ and $-q$ to the first $_4\phi_3$ in the above
equation. The 
result is
\bea\lefteqn{
\kappa(q^{1/2}, b, -c, -q^{1/2})\frac{(-q)^k(1-bc)(q,-bc)_k}{
(1-bcq^{2k})(bc,-q)_k}\,A_k }\\
& & =\sum_{r=0}^{k} \frac{(q^{-k}, bcq^k)_r\; q^r}{
(q,-q,cq^{1/2}, -bq^{1/2})_r}
\sum_{s=0}^{n-1} \frac{(q^{1-n}, bcq^{n+1})_s\,q^s}{(q,-q,
bq^{3/2},-cq^{3/2})_s}
\nonumber \\
& & \mbox{\hspace{0.2in}}\cdot
\int_0^\pi\frac{(e^{2i\t},e^{-2i\t})_\infty\,d\t}{h(\cos\t;
q^{s+1/2}, bq, -cq, -q^{r+1/2})}. \nonumber
\eea
The integral on the right-hand side of the formula (2.5) is
\[
\frac{2\pi(bcq^{r+s+3})_\infty}{(q,bq^{s+3/2},-cq^{s+3/2},-q^{r+s
+1},
-bcq^2, -bq^{r+3/2}, cq^{r+3/2})_\infty},
\]
which can be written as
\[
\kappa(q^{1/2}, b,-c,-q^{1/2})\frac{(1+bc)(1+bcq)(1-qb^2)(1-qc^2)}{
(bcq)_2(bcq^3)_{r+s}}
\]\dis\dis
\[
\cdot (-bq^{3/2}, cq^{3/2})_r\,(bq^{3/2}, -cq^{3/2})_s (-q)_{r+s}.
\]
Therefore (2.5) leads to
\be
\frac{(1-bcq)(1-bcq^2)}{(-bc)_2(1-qb^2)(1-qc^2)}
\frac{(1-bc)(q,-bc)_k}{(1-bcq^{2k})(bc,-q)_k} (-q)^k A_k
\ee
\dis\dis
\[
=\sum_{r=0}^k
\frac{(q^{-k},bcq^k,
-bq^{3/2},cq^{3/2})_r}{(q,bcq^3,cq^{1/2},-bq^{1/2})_r}
q^r{}_3\phi_2\left(\left.\ba{c}
q^{1-n},bcq^{n+1},-q^{r+1} \\
bcq^{r+3}, \qquad -q \ea \right|\, q, q\right).
\]
The above $_3\phi_2$ can be summed by the $q$-analog of the
Pfaff-Saalsch%
\"{u}tz theorem, (II.12) in \cite{Ga:Ra}. It's sum is
\[
\frac{(-bcq^2, q^{2+r-n})_{n-1}}{(-q^{1-n}, bcq^{r+3})_{n-1}}
\]
which clearly  vanishes if $r\le n-2
$. Thus $A_k=0$ if $k\le n-2$. When $k\ge n-1$, replace $r$ by
$r+n-1$ in
(2.6) and simplify the result to see that the right-hand side of
(2.6) is
\[
\frac{(-bcq^2,q^{-k},bcq^k,
-bq^{3/2},cq^{3/2})_{n-1}}{(-q^{1-n},bcq^3,cq^{1/2},-bq^{1/2})_{n
-1}}
\frac{q^{n-1}}{(bcq^{n+2})_{n-1}}
\]
\dis\dis
\[
\cdot{}_4\phi_3\left(\left. 
\ba{c} q^{n-1-k}, bcq^{n-1+k}, -bq^{n+1/2}, cq^{n+1/2} \\
bcq^{2n+1},\quad -bq^{n-1/2},\quad cq^{n-1/2} \ea
\right|\,q,q\right).
\]
When $k=n-1$ the above ${}_4\phi_3$ is 1. If $k=n,n+1$, the
aforementioned
$_4\phi_3$ has only 1, 2 terms; respectively. After some
simplification we
find
\be
A_n=\frac{(1-bq^{1/2})(1+cq^{1/2})(1+bcq^{2n})(1+bcq^n)(1+q^n)(1-
b/c)}{
(1-bcq^{2n-1})(1-bcq^{2n+1})}\,cq^{-1/2},
\ee
\be
A_{n+1}=-\frac{(1-bq^{1/2})(1+cq^{1/2})(1-bq^{n+1/2})(1+cq^{n+1/2
})(1+q^n)(1+q^{n+1})}{
(1-bcq^{2n})(1-bcq^{2n+1})}\,(bcq^{-1}),
\ee
\be
A_{n-1}=\frac{(1-bq^{1/2})(1+cq^{1/2})(1+bq^{n-1/2})
(1-cq^{n-1/2})(1+bcq^{n-1})(1+bcq^n)}{
(1-bcq^{2n-1})(1-bcq^{2n})}.
\ee

Using the dual relationships mentioned at the beginning of this
section we 
see that Theorem 2.1 is equivalent to the following theorem.

\begin{thm}
The connection coefficient formula 
\be
P_n^{(\a, \b)}(x|q) = \frac{q^{-n/2}(1 - q^{\a + \b + n + 1})
(1 - q^{\a + \b + n + 2})}{(- q^{(\a + \b  + 1)/2}; q^{1/2})_2(1 -
q^{n + (\a + \b + 1)/2})(1 - q^{n + (\a + \b + 2)/2})} P_n^{(\a +
1, \b + 1)}(x|q)
\ee
\bea
+ \frac{q^{(\a + \b + 2 -n)/2}(1 - q^{\a + \b + n + 1})
(1 + q^{n + (\a + \b  + 1)/2}) (1 - q^{(\a - \b)/2})}
{(- q^{(\a + \b  + 1)/2}; q^{1/2})_2(1 - q^{n + (\a + \b)/2})(1 -
q^{n + (\a + \b + 2)/2})} P_{n - 1}^{(\a + 1, \b + 1)}(x|q)
\nonumber 
\eea
\bea
- \frac{q^{(3\a + \b + 4 -n)/2}(1 - q^{\a + n})(1 - q^{\b + n})}
{(- q^{(\a + \b  + 1)/2}; q^{1/2})_2(1 - q^{n + (\a + \b)/2})(1 -
q^{n + (\a + \b + 1)/2})} P_{n - 2}^{(\a + 1, \b + 1)}(x|q),
\nonumber 
\eea
holds.
\end{thm}

\newpage
\section{An Eigenvalue Problem}
\setcounter{equation}{0}
In this section we characterize the eigenvalue and the
eigenfunction of
 the eigenvalue problem
\be
(T_{\a,\b;q} g)(x)=\l g(x).
\ee
This characterization is stated as Theorem 3.1 at the end of the
present section.

It is tacitly assumed in (3.1) that $g$ belongs to the domain of
$T_{\a,\b;q}$ and $\l g$ belongs to its range. Now 
assume
\be
g(x):=g(x;\l|q)\sim \Sum a_n(\l|q)P_{n}^{(\a,\b)}(x|q) .
\ee
Since $g\in L_2[w_{\alpha,\beta}(x;q)]$ then (1.10) implies
\be
\sum_{n=0}^{\infty} \, h_n^{(\alpha,\beta)}(q) |a_n(\lambda|q)|^2
\;<\; \infty.
\ee
The condition (3.3) and the eigenvalue equation (3.1) will
characterize the 
eigenvalues $\lambda$ and the eigenfunctions $g$.

It is clear that
\bea\lefteqn{
\l \Sum a_n(\l|q)P_{n}^{(\a,\b)}(x|q)}\\
& &  =\Sum a_n(\l|q)\Int
w_{\alpha + 1, \beta + 1}(x;q)\;
K^{(\a,\b)}(x,t)P_{n}^{(\a,\b)}(t|q)\;dt
  \nonumber \\
&
&=\frac{(q,q^{\a+2},q^{\b+2},-q^{\frac{\a+\b+1}2},-q^{\frac{\a+\b
+2}2)})_\infty
(1-q)}{4\pi
(q^{\frac{\a+\b+4}2},q^{\frac{\a+\b+5}2})_\infty(1-q^{\a+\b+3})}
q^{-\frac14(2\a+1)}\nonumber \\
& & \quad\cdot \Sum
a_n(\l|q)\sum_{k=0}^\infty\frac{1-q^{\a+\b+3+2k}}
{1-q^{\a+\b+2+k}}\,\frac{(q^{\a+\b+3},q,-q^\frac{\a+\b+3}2)_k}{
(q^{\a+2},q^{\b+2},-q^{\frac{\a+\b+5}2})_k}
q^{-\frac{k}2(2\a+1)}P_{n+1}^{(\a,\b)}(x|q) \nonumber \\
& & \cdot\Int w_{\alpha + 1, \beta +
1}(x;q)P_{k}^{(\a+1,\b+1)}(t|q) P_{n}^{(\a,\b)}(t|q)\;dt .\nonumber
\eea
The next step is to evaluate the integral on the extreme right-hand
side of 
(3.4). This will lead to a three term recurrence relation satisfied
 by the $a_n$'s. By (2.3), the integral on the right, denoted
$I_{k,n}$ with $k$ replaced
by
$k-1$, is
\bea\lefteqn{
I_{k,n}=\Int 
w(t;q^{\a/2+1/4},q^{\a/2+3/4},-q^{\b/2+1/4},-q^{\b/2+3/4})}  \\
& & \cdot(1+q^{\frac{\a+\b+1}2})(1+q^{\frac{\a+\b+2}2})\left\{
\frac{(1-q^{\a+k})(1-q^{\b+k})}{(1-q^{\frac{\a+\b}2+k})(1-q^{\frac{\a+\b+1}2+k})
}q^{\frac{k-1}2}P_{k-1}^{(\a,\b)}(t|q)\right. \nonumber \\
& &
+\;\frac{(1+q^{\frac{\a+\b+1}2+k})(1-q^k)(1-q^{\frac{\a-\b}2})}{
(1-q^{\frac{\a+\b}2+k})(1-q^{\frac{\a+\b+2}2+k})}
q^{\frac{\b-\a+k-1}2}P_{k}^{(\a,\b)}(t|q) \nonumber \\
& & \left.-\:\frac{(1-q^k)(1-q^{k+1})}
{(1-q^{\frac{\a+\b+1}2+k})(1-q^{\frac{\a+\b+2}2+k})}
q^{\frac{\b-\a+k-1}2}P_{k+1}^{(\a,\b)}(t|q)\right\}P_{n}^{(\a,\b)
}(t|q)dt. \nonumber 
\eea
Now the orthogonality relation (1.10) is equivalent to
\bea\lefteqn{
\Int 
w(t;q^{\a/2+1/4},q^{\a/2+3/4},-q^{\b/2+1/4},-q^{\b/2+3/4})
P_{n}^{(\a,\b)}(t|q) P_{j}^{(\a,\b)}(t|q)dt } \nonumber  \\
& & =\frac{2\pi
(q^{\frac{\a+\b+2}2},q^{\frac{\a+\b+3}2};q)_\infty}
{(q,q^{\a+1},q^{\b+1},-q^{\frac{\a+\b+1}2},
-q^{\frac{\a+\b+2}2};q)_\infty
}
\nonumber \\
& & \frac{
(1-q^{\a+\b+1})(q^{\a+1},q^{\b+1},-q^{\frac{\a+\b+3}2};q)_j}{
(1-q^{\a+\b+1+2j})(q^{\a+\b+1},q,-q^{(\a+\b+1)/2};q)_j}
q^{\frac{j}{2}(2\a+1)} \delta_{n,j}. \nonumber
\eea
Substituting this in (3.5) we find that
\bea\lefteqn{
I_{k,n}=\frac{2\pi
(q^{\frac{\a+\b+2}2},q^{\frac{\a+\b+3}2};q)_\infty}
{(q,q^{\a+1},q^{\b+1},-q^{\frac{\a+\b+1}2},
-q^{\frac{\a+\b+2}2};q)_\infty
}
(1+q^{\frac{\a+\b+1}2})(1+q^{\frac{\a+\b+2}2})}\\
& & \quad\cdot\frac{(1-q^{\a+\b+1})}{(1-q^{\a+\b+1+2n})}
\frac{
(q^{\a+1},q^{\b+1},-q^{\frac{\a+\b+3}2};q)_n}{(q^{\a+\b+1},q,-q^{
\frac{\a+\b+1}2};q)_n}q^{\frac{n}2(2\a+1)}
\nonumber \\
& &
\quad\cdot\left\{\frac{(1-q^{\a+k})(1-q^{\b+k})}{(1-q^{\frac{\a+\b}2+k})
(1-q^{\frac{\a+\b+1}2+k})}q^{\frac{k-1}2}\delta_{n,k-1}\right.\nonumber\\
& & +\frac{(1+q^{\frac{\a+\b+1}2+k})(1-q^k)(1-q^{\frac{\a-\b}2})}{
(1-q^{\frac{\a+\b}2+k})(1-q^{\frac{\a+\b+2}2+k})}q^{\frac{\b-\a+k
-150z49z}2}\delta_{n,k}
\nonumber \\ & &
-\left.\frac{(1-q^k)(1-q^{k+1})q^{\frac{\b-\a+k-1}2}}{
(1-q^{\frac{\a+\b+1}2+k})(1-q^{\frac{\a+\b+2}2+k})}\delta_{n,k+1}
\right\}.\nonumber
\eea
From (3.4) and (3.6) we have
\bea\lefteqn{
\l \Sum a_n(\l|q)P_{n}^{(\a,\b)}(x|q)}\\
& &= \frac{(1-q)(1-q^{\frac{\a+\b+2}2})(1-q^{\frac{\a+\b+3}2})}{
2(1-q^{\a+1})(1-q^{\b+1})(1-q^{\a+\b+3})}q^{-\frac14(2\a+1)}
(1+q^{\frac{\a+\b+1}2})(1+q^{\frac{\a+\b+2}2})\nonumber \\
& & \quad\cdot \sum_{k = 1}^\infty
\frac{1-q^{\a+\b+1+2k}}{1-q^{\a+\b+1+k}}\frac{
(q^{\a+\b+3},q,-q^{\frac{\a+\b+3}2};q)_{k-1}}{
(q^{\a+2},q^{\b+2},-q^{\frac{\a+\b+5}2};q)_{k-1}}q^{-\frac{k-1}2(
2\a+1)}
P_{k}^{(\a,\b)}(x|q)\nonumber \\
& & \quad\cdot\left\{
\frac{(1-q^{\a+k})(1-q^{\b+k})(q^{\a+1},q^{\b+1},-q^{\frac{\a+\b+
3}2};
q)_{k-1}(1-q^{\a+\b+1})}{(1-q^{\frac{\a+\b}2+k})(1-q^{\frac{\a+\b
+1}2+k})(
q^{\a+\b+1},q,-q^{\frac{\a+\b+1}2};q)_{k-1}(1-q^{\a+\b+2k-1})}\right.
\nonumber \\
& & \quad\cdot q^{\frac{k-1}2+\frac{k-1}2(2\a+1)}a_{k-1}(\l|q)
\nonumber \\
& & +\frac{(1+q^{\frac{\a+\b+1}2+k})(1-q^k)(1-q^{\frac{\a-\b}2})
(q^{\a+1},q^{\b+1},-q^{\frac{\a+\b+3}2};
q)_{k}}{(1-q^{\frac{\a+\b}2+k})(1-q^{\frac{\a+\b+2}2+k})(
q^{\a+\b+1},q,-q^{\frac{\a+\b+1}2};q)_{k}}q^{\frac{k}2(2\a+1)+\frac{\b-\a+k-1}2}
\nonumber \\
& &
\quad\cdot\frac{(1-q^{\a+\b+1})}{(1-q^{\a+\b+1+2k})}\,a_k(\l|q)\nonumber \\
& &
-\frac{(1-q^k)(1-q^{k+1})(q^{\a+1},q^{\b+1},-q^{\frac{\a+\b+3}2};
q)_{k+1}(1-q^{\a+\b+1})}{(1-q^{\frac{\a+\b+1}2+k})(1-q^{\frac{\a+
\b+2}2+k})(
q^{\a+\b+1},q,-q^{\frac{\a+\b+1}2};q)_{k+1}(1-q^{\a+\b+2k+3})}\nonumber \\
& &\quad\cdot \left.\makebox[0in]{\rule{0in}{0.3in}}
q^{\frac{k+1}2(2\a+1)+\frac{\b-\a+k-1}2}\,a_{k+1}(\l|q)\right\}\nonumber 
%
\eea
After some simplification we find
\bea\lefteqn{
\l \sum_{n=1}^\infty a_n(\l|q) P_{n}^{(\a,\b)}(x|q)
 =\sum_{k=1}^\infty
P_{k}^{(\a,\b)}(x|q)\left[\frac{(1-q)(1-q^{\a+\b+k})
}{
2(1-q^{\frac{\a+\b}2+k})(1-q^{\frac{\a+\b-1}2+k})}q^{\frac{k-1}2}
\,a_{k-1}(\l|q)\right.} \\
& & +\frac{(1-q)(1-q^{\frac{\a-\b}2})(1+q^{\frac{\a+\b+1}2+k})}{
2(1-q^{\frac{\a+\b}2+k})(1-q^{\frac{\a+\b+2}2+k})}q^{\frac{\a+\b+
k}2}
\,a_{k}(\l|q)\nonumber \\
& &
-\left.\frac{(1-q)(1-q^{\a+k+1})(1-q^{\b+k+1})q^{(3\a+\b+k+1)/2}}{
2(1-q^{\a+\b+1+k})(1-q^{\frac{\a+\b+2}2+k})(1-q^{\frac{\a+\b+3}2+
k})}
\,a_{k+1}(\l|q)\right]\nonumber
\eea
By equating the coefficients of $P_n^{(\a,\b)}(x|q)$ on both sides
of
(3.8) we establish the following
three-term recurrence relation for $a_n$'s 
\bea\lefteqn{
-\l a_{k}(\l|q) q^{\frac{\a}2+\frac14} =
\frac{(1-q)(1-q^{\a+k+1})(1-q^{\b+k+1})\; q^{(3\a + \b + k +
1)/2}}{
2(1-q^{\a+\b+1+k})(1-q^{\frac{\a+\b+2}2+k})(1-q^{\frac{\a+\b+3}2+
k})}
\,a_{k+1}(\l|q)}\\
& & -\frac{(1-q)(1-q^{\frac{\a-\b}2})(1+q^{\frac{\a+\b+1}2+k})}{
2(1-q^{\frac{\a+\b}2+k})(1-q^{\frac{\a+\b+2}2+k})}q^{\frac{\a+\b+
k}2}
\,a_{k}(\l|q)\nonumber \\
& & -\frac{(1-q)(1-q^{\a+\b+k})
}{
2(1-q^{\frac{\a+\b}2+k})(1-q^{\frac{\a+\b-1}2+k})}q^{\frac{k-1}2}
\,a_{k-1}(\l|q), \quad k>0.\nonumber
\eea
The  limiting case $q\to 1^-$  of the recursion relation (3.9) is
\bea\lefteqn{
-\l
a_{k}(\l)=\frac{2(\a+1+k)(\b+1+k)}{(\a+\b+1+k)(\a+\b+2+2k)_2}a_{k
+1}(\l)
}\\
& & +\frac{2(\b-\a)}{(\a+\b+2k)(\a+\b+2+2k)} a_{k}(\l)
-\frac{2(\a+\b+k)}{(\a+\b+2k-1)(\a+\b+2k)} a_{k-1}(\l).
\nonumber
\eea
for $k>0$, which is (4.16) of \cite{Is:Zh}

It is clear from (3.3) that $a_0(\l|q)=0$ and that $a_1(\l|q)$ is
arbitrary. It is also clear from (3.3) that $a_k(\l|q)/a_1(\l|q)$
is a
polynomial in $\l$ of degree $k-1$. It is more convenient to
renormalize
$a_k(\l|q)/a_1(\l|q)$ in terms of monic polynomials. Thus we set
\be
a_{k+1}(\l|q)=\frac{(q^{\a+\b+2},q^{\frac{\a+\b+4}2},q^{\frac{\a+
\b+5}2};q)_k}{
(q^{\a+2},q^{\b+2};q)_k}(-1)^kb_k(\frac{2\l
q^{1/2}}{1-q})q^{-({k^2}/4+(
\a+\b/2+1)k)}.
\ee
There is no loss of generality in taking $b_0(2\l/(1-q))=1$. In
terms of
the $b_n$'s, (3.8) becomes
\bea\lefteqn{
{b_{k+1}(\mu)}={b_{k}(\mu)}\left[\mu+\frac{
(1-q^{\frac{\b-\a}2})(1+q^{\frac{\a+\b+3}2+k})}{(1-q^{\frac{\a+\b
+2}2+k})
(1-q^{\frac{\a+\b+4}2+k})}q^{\a/2+3/4+k/2}\right]}\\
& & +\frac{(1-q^{\a+1+k})(1-q^{\b+1+k})q^{k+\frac{\a+\b}2+1}}{
(1-q^{\frac{\a+\b+1}2+k})(1-q^{\frac{\a+\b+2}2+k})^2(1-q^{\frac{\a+\b+3}2+k})}
{b_{k-1}(\mu)},\nonumber
\eea
where 
\be
\mu=\frac{2\l q^{1/2}}{1-q},\quad b_{-1}(\mu)=0,\quad b_0(\mu)=1.
\ee

In Section 5 we shall determine the large $n$ behavior of the
polynomials $b_n(x)$
and $a_n(\lambda|q)$.  These asymptotic results will
 be used to prove the following theorem.
\begin{thm}
The  eigenvalue problem (3.1)-(3.3)  has a countable infinite
number 
of eigenvalues. The eigenvalues are  $(1-q)/2$ times the 
reciprocals of the 
roots of the transcendental equation
\be
(-p^{\a +3/2}(1-q)x/2;p)_{\infty}{} _2\phi_1\left(\left.\ba{c}
p^{\a +1},(1-q)x p^{1/2}/2\\
- p^{\a+ 3/2}(1-q)x/2\ea\,\right|\;p,p^{\b+1}\right) = 0,\; \; p :=
q^{1/2}.
\ee
Furthermore $\lambda = 0$ is not an eigenvalue and the eigenspaces
are
one dimensional.\end{thm}
\begin{thm}
 Any eigenfunction  $g(x; \lambda|q)$ corresponding to an
eigenvalue $\lambda$ 
is a constant multiple of ${\cal E}_q(x;-i,\l)$.
\end{thm}

\newpage
\section{A q-Analog of Wimp's Polynomials.}
\setcounter{equation}{0}
In this section we find an explicit solution of (3.12).
\begin{thm}
The polynomial $\{b_n(x)\}$ generated by (3.12) and (3.13) are
given by
\bea\lefteqn{
b_n(\mu)=\sum_{j=0}^n\frac{(p^{-\b-n-1},-p^{-\a-n-1};p)_j}{
(p,p^{-2n-\a-\b-2};p)_j}(-1)^jp^{j/2}\mu^{n-j}}\\
& & \mbox{\hspace{0.7in}}\cdot{}_4\phi_3\left(\left.\ba{c}
p^{-j},p^{2n+\a+\b+3-j},p^{\b+1},-p^{\a+1}\\
p^{\a+\b+2},p^{n+\b+2-j},-p^{\a+n+2-j}\ea\,\right|\;p,p\right),
\nonumber
\eea
where
\be
p:=q^{1/2}.
\ee
\end{thm}
{\bf Proof.} From (4.1) it is clear that $b_0(\mu) = 1$ and that
$b_1(\mu)$
satisfies (3.12) when $k = 0$. Since the solution of the initial
value 
problem (3.12)-(3.13) must be unique, all we need to do is to
verify that
the right-hand side of (4.1) satisfies (3.12). The actual process
of this
 verification is rather long and tedious.

First, we shall rewrite (4.1) in the form 
\be
b_n(\mu)=\sum_{j=0}^n\sum_{k=0}^j A_k\, B_{j-k}^{(n)}\,
(-1)^{j+k}p^{j/2}\mu^{n-j},
\ee
where
\be
A_k=\frac{(p^{\b+1},-p^{\a+1};p)_k}{(p,p^{\a+\b+2};p)_k}, \quad
B_k^{(n)} = 
\frac{(p^{-\b-n-1},-p^{-\a-n-1};p)_{k}}{
(p,p^{-2n-\a-\b-2};p)_{k}}.
\ee
Since $A_0 = B_0^{(n)} = 1$, we then have
\begin{eqnarray}
b_{n+1}(\mu) =  \mu^{n+1}+\sum_{j=1}^{n+1}\sum_{k=0}^j A_k\,
B_{j-k}^{(n+1)}\, 
(-1)^{j+k}p^{j/2}\mu^{n+1-j} 
\eea
\bea
= \mu\left[
b_n(\mu)-\sum_{j=1}^n\sum_{k=0}^j A_k\,
B_{j-k}^{(n)}\,(-1)^{j+k}p^{j/2}\mu^{n-j}
+\sum_{j=0}^{n+1}\sum_{k=0}^j A_k\, B_{j-k}^{(n+1)}\,
(-1)^{j+k}p^{j/2}\mu^{n+1-j}\right], \nonumber
\eea
where the last line is obtained by separating $\mu^{n}$ from the 
rest of the series on the right-hand side of (4.1). Our first 
aim is to bring the same factor of $b_n(\mu)$ as is shown in
(3.12), so we 
make a further separation  of the series on (4.5) and find that
\be
b_{n+1}(\mu) = b_n(\mu)\left[\mu +
\frac{(1-p^{\b-\a})(1+p^{\a+\b+3+2n})}
{(1-p^{\a+\b+2+2n})(1-p^{\a+\b+4+2n})}p^{\a+n+3/2}\right]\; +
c_n(\mu),
\ee
where
\bea \lefteqn{
c_n(\mu) = \sum_{j=0}^{n-1}\sum_{k=0}^{j+1}
A_k\, B_{j+1-k}^{(n)}(-1)^{j+k}p^{\frac{j+1}2} 
\mu^{n-j} 
-\sum_{j=0}^{n}\sum_{k=0}^{j+1}
A_kB^{(n+1)}_{j+1-k}(-1)^{j+k}p^{\frac{j+1}2}
\mu^{n-j} } \\
& & -\frac{(1-p^{\b-\a})(1+p^{\a+\b+3+2n})p^{\a+n+1}}
{(1-p^{\a+\b+2+2n})(1-p^{\a+\b+2n+4})}\sum_{j=0}^{n}\sum_{k=0}^{j}
A_kB_{j-k}^{(n)}(-1)^{j+k}p^{\frac{j+1}2}
\mu^{n-j}. \nonumber
\eea
The rest of the exercise is to show that $c_n(\mu)$ is actually a
multiple
of $b_{n-1}(\mu)$, the same multiple as in (3.12). The coefficient
of 
$\mu^{n}$ in $c_n(\mu)$ is
\bea
p^{j/2}\sum_{k=0}^{1}[B_{1-k}^{(n)} - B_{1-k}^{(n+1)}](-1)^kA_k 
-\frac{(1-p^{\b-\a})(1+p^{\a+\b+3+2n})p^{\a+n+3/2}}
{(1-p^{\a+\b+2+2n})(1-p^{\a+\b+2n+4})}, \nonumber
\eea
which vanishes by the use of (4.4), verifying that $c_n(\mu)$ is 
a polynomial of degree $n-1$ in $\mu$. We thus have
\bea
c_n(\mu) = \frac{(1+p^{\a+\b+3+2n})(1-p^{\b-\a})p^{\a+n+1}}
{(1-p^{\a+\b+2+2n})
(1+p^{\a+\b+4+2n})} 
\sum_{j = 0}^{n-1}\sum_{k=0}^{j+1}(-1)^{j+k} A_k B_{j+1-k}^{(n)}
p^{(j+2)/2} \mu^{n-j-1}
+ d_n(\mu),
\eea
where
\be
d_n(\mu):= \sum_{j=0}^{n-2} \sum_{k=0}^{j+2} (-1)^{j+k+1} A_k
B_{j+2-k}^{(n)}
p^{(j+2)/2} \mu^{n-j-1}
\ee
\bea
\qquad \qquad - \sum_{j=0}^{n-1} \sum_{k=0}^{j+2} (-1)^{j+k+1} A_k
B_{j+2-k}^{(n+1)}
p^{(j+2)/2} \mu^{n-j-1}. \nonumber
\eea
Since
\bea \lefteqn{
B_{j+2-k}^{(n)} -  B_{j+2-k}^{(n+1)} = \frac{(p^{-\b-n-1},-
p^{-\a-n-1};p)_{j+1-k}}
{(p;p)_{j+1-k}(p^{-2n-\a-\b-4};p)_{j+4-k}}} \\
& & \cdot \left[p^{-2n-\a-\b-3}(1-p^{j-k+1})(1 - p^{-2n-\a -\b -4})
-p^{-n-\b-2} (1-p^{\b-\a})(1 -
p^{-4n-2\a-2\b-5+j-k})\right],\nonumber
\eea
we find that the coefficient of $\mu^{n-1-j}$ in $d_n(\mu)$ is, for

$0 \le j \le n-2$,
\bea \lefteqn{
\frac{(1-p^{-\b-n-1})(1+p^{-\a-n-1})p^{-2n-\a-\b-2}}{(p^{-2n-\a-\b-3};p)_3}
\sum_{k=0}^{j}A_k B_{j-k}^{(n-1)} (-1)^{j+k}
p^{j/2}-(1-p^{\b-\a})p^{-\b-n-2}}
\\
& & \cdot \sum_{k=0}^{j+1} A_k(1-p^{-4n-2\a-2\b-5+j-k}) (-1)^{j+k}
p^{(j+2)/2}
\frac{(p^{-\b-n-1}, -p^{-\a-n-1};p)_{j+1-k}}{(p;p)_{j+1-k}
(p^{-2n-\a-\b-4)};p)_{j+4-k}}.\nonumber
\eea

Now we combine the second series in (4.11) with the coefficients of
$\mu^{n-1-j}$ in the  first  series on the right-hand side of (4.8)
which, by virtue of the identity
\bea\lefteqn{
\frac{(1+p^{\a+\b+3+2n})p^{\a+n+1}}{
(1-p^{\a+\b+2+2n})(1-p^{\a+\b+4+2n})(p^{-2n-\a-\b-2};p)_{j+1-k}}
- \frac{p^{-\b-n-2}(1-p^{-4n-2\a-2\b-5+j-k})}{
(p^{-2n-\a-\b-4};p)_{j+4-k}}    }\nonumber \\
& & = - \frac{(1-p^{j+1-k})p^{-\b-n-2}}{(1 - p^{\a + \b +
2+2n})(p^{-2n-\a-\b-3};p)_{j+3-k}},\nonumber
\eea
results in the series
\bea
-\frac{(1-p^{\b-\a})(1 - p^{-\b-n -1})(1+p^{-\a -n -
1})p^{-\b-n-1}}
{(1 - p^{2n+\a+\b+2})(p^{-2n-\a -\b - 3};p)_3}
\sum_{k=0}^{j} A_k \, B_{j-k}^{(n-1)}\, (-1)^{k+j} p^{j/2}.
\nonumber
\eea
Adding this to the first series in (4.11) we find, after a 
straightforward calculation, that the coefficient of $\mu^{n-j-1}$,
 $0 \le j \le n-2$ in $c_n(\mu)$ of (4.8) is
\bea
\frac{(1-p^{2\a + 2n + 2})(1-p^{2\b + 2n + 2})p^{\a +\b + 2n + 2}}
{(1-p^{2n+ \a + \b + 3})(1-p^{2n+ \a + \b + 1})(1-p^{2n+ \a + \b +
2})^2}
\sum_{k=0}^{j}A_k\, B_{j-k}^{(n-1)} \,(-1)^{k+j} \, p^{j/2}.
\eea
Finally, collecting the $j = n-1$ term from the series in (4.8) and
(4.9) we find that the constant term in (4.8) is given by
\bea\lefteqn{
\frac{(1-p^{\b -\a})(1+p^{\a +\b + 2n + 3})p^{\a +n + 1}}
{(1-p^{2n+ \a + \b + 2})(1-p^{2n+ \a + \b + 4})}
\sum_{k=0}^{n}A_k\, B_{n-k}^{(n)} \,(-1)^{n+k-1} \, p^{(n+1)/2}}\\
& & - \sum_{k=0}^{n+1}A_k\, B_{n+1-k}^{(n+1)} \,(-1)^{n+k} \,
p^{(n+1)/2}
=: f_n, \nonumber
\eea
say. By (4.4), we get
\bea\lefteqn{
f_n = \frac{(p^{-\b-n-1},-p^{-\a-n-1};p)_{n}
}{(p;p)_{n+1}(p^{-2n-\a-\b-4};p)_{n+2}}(-1)^{n-1}p^{n+1)/2} }\\
& & \cdot [(1-p^{-\b-n-2})(1+p^{-\a -n
-2})(1-p^{-n-\a-\b-3})\nonumber\\
& & \cdot {}_4\phi_3\left(\left.\ba{c}
p^{-n-1},p^{n+\a+\b+4},p^{\b+1},-p^{\a+1} \nonumber \\
p^{\a+\b+2},\qquad p^{\b+2},\qquad -p^{\a+2}\ea\, \right| 
p,p\right)  \nonumber
\\ 
& & +
\frac{(1-p^{\b-\a})(1-p^{2\a+2\b+6+4n})(1-p^{n+1})p^{-\a-2\b-3n-6
}}{(1-p^{\a+\b+2+2n})} \nonumber
\\
& &\mbox{\hspace{0.3in}}{}_4\phi_3\left(\left.\ba{c}
p^{-n},p^{n+\a+\b+3},p^{\b+1},-p^{\a+1} \nonumber \\
p^{\a+\b+2},\qquad p^{\b+2},\qquad -p^{\a+2}\ea\, \right| 
p,p\right) ].
\nonumber 
\eea
To simplify the expression on the right side of (4.14) we first
denote
\[
\phi_n := {}_4\phi_3\left(\left. \ba{c}
p^{-n},p^{n+\a+\b+3},p^{\b+1},-p^{\a+1}\\
p^{\a+\b+2},\; p^{\b+2},\; -p^{\a+2}\ea\, \right|\: p,p\right)
\]
and then use the contiguous relation of Askey and Wilson, see
\cite[Ex 7.5]{Ga:Ra}:
\be
\phi_{n+1} = - \frac{B}{A}\phi_n - \frac{C}{A}\phi_{n-1},
\ee
where
\bea\lefteqn{ 
A=p^{\a+\b-3n+4}(1-p^{n+\a+\b+3})(1-p^{2n+\a+\b+2})(1-p^{n+\a+\b+
2})}
\\
& & \cdot (1-p^{n+\b+2})(1+p^{n+\a+2}), \nonumber
\eea
\bea
C=-p^{2\a+2\b+6-3n}(1-p^n)(1-p^{2n+\a+\b+4})(1-p^{n+1})(1-p^{n+\a
+1})
(1+ p^{n+\b+1}),\nonumber
\eea
\bea
B=-C-A+p^{\a+\b-3n+4}(p^{2n+\a+\b+2})_3(1-p^{\b+1})(1+p^{\a+1}).
\nonumber
\eea
Substituting (4.15) in (4.14)we find after some simplification that
the
coefficients of $\phi_n$
cancel out, so that 
\bea
f_n = -\frac{C}{A}
(-1)^{n-1}p^{(n+1)/2}\frac{(p^{-\b-n-2},-p^{-\a-n-2};p)_{n+1}
}{(p,p^{-2n-\a-\b-4};p)_{n+1}}\,  \phi_{n-1}
\eea
\bea
 = \frac{(1-p^{2\a + 2n + 2})(1-p^{2\b + 2n + 2})p^{\a +\b + 2n +
2}}
{(1-p^{2n+ \a + \b + 3})(1-p^{2n+ \a + \b + 1})(1-p^{2n+ \a + \b +
2})^2}
\sum_{k=0}^{n-1}A_k\, B_{n-1-k}^{(n-1)} \,(-1)^{k+n-1} \,
p^{(n-1)/2}.
\nonumber 
\eea
Combining (4.12) and (4.17) we find that $c_n(\mu)$ is the same as 
the second term on the right side of (3.12) (with $k$ replaced by
$n$).
 This completes the proof of Theorem 4.1.

\newpage
\section{Properties of $\{b_n(x)\}$.}
\setcounter{equation}{0}
In this section we derive asymptotic formulas for the polynomials 
$\{b_n(x)\}$ in different parts of the complex $x$-plane.
We also investigate a closely related set of orthogonal polynomials
and
record their associated continued $J$-fraction.

Recall that we normalized $\{b_n(x)\}$ of (4.1) by 
\be
b_0(x):=1.
\ee
To exhibit the dependence of $b_n(x)$ on the parameters $\a$ and
$\b$
we shall use the  notation $b_n^{(\a,\b)}(x)$ instead of $b_n(x)$.
Our first result concerns the limiting behavior of
$\{b_n^{(\a,\b)}(x)\}$.
\begin{thm}
The limiting relation
\be
\lim_{n\to\infty} x^{n}\,b_n^{(\a,\b)}(1/x)=\frac{(p^{\b+1},
-p^{\a+ 3/2}x;p)_\infty}
{(p^{\a+\b + 2}; p)_{\infty}}
{}_2\phi_1\left(\left. \ba{c} p^{\a+1}, p^{1/2}x \\ -p^{\a+3/2}x
\ea \right|p,p^{\b+1}\right),
\ee
holds uniformly on compact subsets of the complex $x$-plane.
\end{thm}

{\bf Proof.} We may use Tannery's theorem (the discrete version of
the
Lebesgue bounded convergence theorem) to let $n\to\infty$ in (4.1)
after
multiplying it by $x^{-n}$. Thus the sequence
$\{x^{-n}b_n^{\a,\b}(x)\}$ will 
have a finite limit if the series
\[
\sum_{j=0}^\infty \frac{(-x)^{-j} p^{j^2/2}}{(p;p)_j}
{}_3\phi_2\left(\left. \ba{c} p^{-j},p^{\b+1}, -p^{\a+1} \\ 0,\quad
p^{\a+\b+2} 
\ea \right|p,p\right)
\]
converges. Therefore
\be
\lim_{n\to\infty} x^{-n} b_n^{(\a,\b)}(x)=\sum_{j=0}^\infty
\frac{(-x)^{-j} p^{j^2/2}}{(p;p)_j}
\sum_{k=0}^j \frac{ (p^{-j},p^{\b+1}, -p^{\a+1};p)_k}{ (p,\quad
p^{\a+\b+2} 
;p)_k}\,p^k,
\ee
if the right-hand side  exists. In the above sum interchange the
$j$ and
$k$ sums and replace $j$ by $j+k$ to see that the right-hand side
of (5.3)
is
\[
\sum_{k=0}^\infty \frac{(p^{\b+1}, -p^{\a+1};p)_k}{ (p,\; p^{\a
+\b+2}
;p)_k}x^{-k}p^{k/2}\,\sum_{j=0}^\infty\frac{(-x)^{-j}}{(p;p)_j}p^
{j^2/2}.
\]
The $j$ sum is $(p^{1/2}/x;p)_\infty$ by Euler's sum
\cite[(II.2)]{Ga:Ra}.
This shows that 
\be
\lim_{n\to\infty} x^{-n}\,b_n^{(\a,\b)}(x)=(p^{1/2}/x;p)_\infty
{}_2\phi_1\left(\left. \ba{c} -p^{\a+1}, p^{\b+1} \\ p^{\a+\b+2} 
\ea \right|p,p^{1/2}/x\right),
\ee
uniformly on compact subsets of 
the open  disc $\{x: |x| < p^{1/2}\}$. On the other hand from
Theorem
8.1 we know that 
$$\lim_{n\to\infty} x^{n}\,b_n^{(\a,\b)}(1/x)$$
exists  uniformly on compact subsets of the complex plane 
and is an entire function of $x$. The Heine transformation,
\cite[(III.1)]{Ga:Ra}
\be
{}_2\phi_1(a,b;c;q,z)=\frac{(b,az;q)_\infty}{(c,z,q)_\infty}
{}_2\phi_1(c/b,z;az;q,b),
\ee
implies
$$(p^{1/2}/x;p)_\infty
{}_2\phi_1\left(\left. \ba{c} -p^{\a+1}, p^{\b+1} \\ p^{\a+\b+2} 
\ea \right|p,p^{1/2}/x\right)
=\frac{(p^{\b+1}, -p^{\a+ 3/2}/x;p)_\infty}
{(p^{\a+\b + 2}; p)_{\infty}}
{}_2\phi_1\left(\left. \ba{c} p^{\a+1}, p^{1/2}/x \\ -p^{\a+3/2}/x
\ea \right|p,p^{\b+1}\right).$$
Therefore (5.2) holds in the interior of  $\{x: |x|=p^{1/2}\}$ and
analytic
continuation establishes  the validity of (5.2) on compact subsets
of the 
complex plane. This  completes the proof.

We next determine the asymptotic behavior of $b_n^{(\a,\b)}(x)$ at
$x=0$
and in $\{x:\; 0<|x|\le p^{1/2}\}$.
\begin{thm}
We have for $\a\ne \b$,
\be
b_n^{(\a,\b)}(0)\approx C\,p^{n^2/2}\, u^n\qquad \mbox{ as
}n\to\infty,
\ee
where $C$ is a nonzero constant and $|u|<1$.
\end{thm}
{\bf Proof.} Clearly
\bea\lefteqn{
b_n^{(\a,\b)}(0)=\frac{(p^{-\b-n-1},
-p^{-\a-n-1};p)_n}{(p,p^{-2n-\a-\b-2};p)_n}\,(-1)^np^{n/2}
} \\
& & \cdot{}_4\phi_3\left(\left. \ba{c}
p^{-n}, p^{n+\a+\b+3},p^{\b+1},-p^{\a+1} \\
p^{\a+\b+2}, p^{\b+2},-p^{\a+2} \ea \right|p,p\right).
\nonumber
\eea
Ismail and Wilson \cite{Is:Wi} proved that if $|z|<1$ then
\be
{}_4\phi_3\left(\left. \ba{c}
q^{-n},abcdq^{n-1},az, a/z \\
ab,ac,ad \ea \right|q,q\right)\approx
\left(\frac{a}z\right)^n\frac{
(az,bz,cz,dz;q)_\infty}{(z^2, ab,ac,ad;q)_\infty},
\ee
as $n\to\infty$. Now apply (5.8) with $a=ip^{1+(\a+\b)/2},\;
b=-ip^{1+(\a+\b)/2}$, $c=-ip^{1-(\b-\a)/2}$,\\ 
$d=ip^{1-(\a-\b)/2}$,
\[
z=\left\{\ba{cl} ip^{(\a-\b)/2} & \mbox{ if } \a>\b \\
-ip^{(\b-\a)/2} & \mbox{ if } \b>\a. \ea \right.
\]
Therefore
\[
b_n^{(\a,\b)}(0)\approx
\frac{(p^{\b+2},-p^{\a+2};p)_n}{(p,p^{n+\a+\b+3};p)_n}
(-1)^np^{n^2/2}\left(\frac{a}z\right)^n\frac{
(az,bz,cz,dz;q)_\infty}{(z^2,p^{\a+\b+2},
p^{\b+2},-p^{\a+2};p)_\infty},
\]
which implies (5.6).
\begin{cor}
The values of $C$ and $u$ in (5.6) are given by
\be
u =\left\{\ba{cl} p^{\b+1} & \mbox{ if } \a>\b \\
-p^{\a+1} & \mbox{ if } \b>\a \ea \right.,\qquad
C=\left\{\ba{cl}
\frac{(-p^{\a+1},p^{\a+1};p)_\infty}{(1+p^{\a-\b})(p^{\a+\b+2};p)
_\infty}
&  \mbox{ if } \a>\b \\
\frac{(-p^{\b+1},p^{\b+1};p)_\infty}{(1+p^{\b-\a})(p^{\a+\b+2};p)
_\infty}
&  \mbox{ if } \b>\a. \ea \right.
\ee
\end{cor}

The  only case left now is to determine the large $n$ behavior of
 $b_n^{(\a,\b)}(0)$
on the zeros of 
\be
F(x):=\frac{(p^{\b+1},-p^{\a+3/2}/x;p)_\infty}{(p^{\a+\b+2};p)_\i
nfty}
\,{}_2\phi_1\left(\left. \ba{c}
p^{\a+1}, p^{1/2}/x \\ -p^{\a+3/2}/x \ea \right|p,p^{\b+1}\right).
\ee
This is not straightforward and requires some preliminary results.

\begin{thm}
The function
\bea\lefteqn{
Y_k^{(\a,\b)}(x)=(-x)^{-k}\frac{(p^{2\a+4},p^{2\a+4};p^2)_k\,p^{k
(\a+\b+3+k)}}{
(p^{\a+\b+3},p^{\a+\b+4},p^{\a+\b+4},p^{\a+\b+5};p^2)_k} }\\
& & 
\cdot {}_2\phi_1\left(\left. \ba{c}
-p^{\a+2+k}, p^{\b+2+k} \\p^{\a+\b+2k+4} \ea
\right|p,\frac{p^{1/2}}x
\right),
\nonumber
\eea
satisfies the three term recurrence relation (3.12), with
$p=q^{1/2}$.
\end{thm}

To prove Theorem 5.4 we used MACSYMA to first
find the multiple of the ${}_2\phi_1$ then verified that the
$Y_k^{(\a, \b)}$ of
(5.11) indeed satisfies (3.12) by equating coefficients of powers
of
$1/x$.  Theorem 5.4 is also a limiting case of a result of Gupta,
Ismail 
and Masson \cite{Gu:Is}, as will be explained in \S 8. 

It readily follows from Theorem 5.4 that
\[
X_\nu^{(\a,\b)}(x)=(-x)^{-\nu}(p^{1/2}/x;p)_\infty
{}_2\phi_1\left(\left. \ba{c}
-p^{\a+2+\nu}, p^{\b+2+\nu} \\p^{\a+\b+2\nu+4} \ea \right|
p,\frac{p^{1/2}}{x}
\right)
\]
satisfies the three term recurrence relation
\bea\lefteqn{
\frac{(1-q^{\a+2+\nu})(1-q^{\b+2+\nu})\,q^{\nu+(\a+\b+4)/2}}{
(1-q^{\nu+(\a+\b+3)/2})(1-q^{\nu+(\a+\b+4)/2})^2(1-q^{\nu+(\a+\b+
5)/2})}
X_{\nu+1}^{(\a,\b)}(x) }\\
& &
=X_\nu^{(\a,\b)}(x)\left[x+\frac{(1-q^{(\b-\a)/2})(1+q^{\nu+(\a+\b+3)/2})}
{(1-q^{\nu+(\a+\b+2)/2})(1-q^{\nu+(\a+\b+4)/2})}\,q^{(\nu+\a+3/2)
/2} \right]
+X_{\nu-1}^{(\a,\b)}(x). \nonumber
\eea
The Heine transformation (5.5) yields the alternate representation
\be
X_{\nu}^{(\a,\b)}(x) = 
(-x)^{-\nu} \frac{(p^{\b +\nu + 2}, -p^{\a +\nu +
5/2}/x;p)_{\infty}}
{(p^{\a + \b + 2\nu + 4};p)_{\infty}} 
\ee
 $$\cdot {}_2\phi_1\left(\left. \ba{c}
p^{\a+ \nu +2}, p^{1/2}/x \\- p^{\a+\nu +5/2}/x \ea \right|p, p^{\b
+ \nu + 2}
\right).$$

 According to Theorem 4.5 of \cite{Is:Zh}
\be
C_{\nu}C_{\nu+1}\cdots
C_{\nu+n-1}X_{\nu+n}^{(\a,\b)}(x)=b_n^{(\a+\nu,\b+\nu)}
(x)
X_\nu^{(\a,\b)}(x)+b_{n-1}^{(\a+\nu+1,\b+\nu+1)}(x)X_{\nu-1}^{(\a
,\b)}(x),
\ee
where $C_{\nu}$ is the coefficient of $X_{\nu+1}(x)$ in (5.12), see
(4.28)
in \cite{Is:Zh}.
\begin{thm}
Let $\xi$ be a zero of $F(x)$ of (5.10). The large $n$ behavior of
$b_n^{(\a,\b)}(\xi)$ is
\be
b_n^{(\a,\b)}(\xi)\approx
\frac{\xi^{-n}(q^{\a+2},q^{\b+2};q)_\infty
q^{n(n+\a+\b+3)/2}}{
(q^{(\a+\b+3)/2},q^{(\a+\b+5)/2};q)_\infty\,(q^{(\a+\b+4)/2};q)_\infty^2}
\;\frac{1}{X_0^{(\a,\b)}(\xi)}.
\ee
\end{thm}
{\bf Proof.} It is clear from (5.10) and (5.13) that $F(x)=0$ if
and only if
$X_{-1}^{(\a,\b)}(x)=0$. The recurrence relation (5.14) shows that
if $X_{\nu}^{(\a, \b)}(x)$ and  $X_{\nu -1}^{(\a, \b)}(x)$ vanish
at $x = \zeta$ then  $X_{\nu +n}^{(\a, \b)}(\zeta) = 0$ for all
$n$, $ n = 2, 3, \dots$. But it is obvious that
\[
X_n^{(\a,\b)}(x)\approx (-x)^{-n}, \qquad as\;  n\to\infty.
\]
Thus  $X_{\nu}^{(\a, \b)}(x)$ and  $X_{\nu -1}^{(\a, \b)}(x)$ have
no common zeros. Now (5.14) implies
\be
b_n^{(\a,\b)}(\xi)\approx \frac{(q^{\a+2},q^{\b+2};q)_{\infty}
q^{n(n+\a+\b+3)/2}X_n^{(\a,\b)}(\xi)}{
(q^{(\a+\b+3)/2},q^{(\a+\b+5)/2};q)_{\infty}\,(q^{(\a+\b+4)/2};q)
_{\infty}^2
X_0^{(\a,\b)}(\xi)},
\ee
 and we have established (5.15).

We are now in a position to prove  Theorem 3.1.

\bigskip
{\bf Proof of Theorem 3.1}. From (1.12) and (3.11) it follows that
\be
h_{n+1}^{(\a,\b)}(q)\,|a_{n+1}(\lambda|q)|^2 = O(q^{-n(n+3\a+2\b
-1/2)/2}
|b_n^{(\a,\b)}(2\lambda q^{1/2}/(1-q))|^2).
\ee
Now (5.2), (5.6) , (5.9) and (5.15) show that
$\sum_1^{\infty} h_{n}^{(\a,\b)}(q)\,|a_{n}(\lambda|q)|^2$
converges  if and only 
if $2q^{1/2}\lambda/(1-q)$ is a zero of $X_{-1}^{(\a,\b)}(x)$. The
eigenspaces
are one dimensional since the eigenfunction with an eigenvalue
$\lambda$ must
be given by
\be
g(x|\lambda) =\sum_1^{\infty} a_{n}(\lambda|q)
P_{n}^{(\a,\b)}(x|q).
\ee

Next we prove that the polynomials $\{i^{-n}b_n^{(\a,\b)}(ix)\}$
are
orthogonal on a bounded countable set when $\alpha$  and $\b$ are
not real
and are complex conjugates.

Set
\be
s_n^{(\a,\b)}(x):=i^{-n}b_n^{(\a,\b)}(ix).
\ee
The three term recurrence relation (3.12) leads to the following
three term
recurrence relation for the $s_n$'s.
\bea\lefteqn{
s_{n+1}^{(\a,\b)}(x)=\left[x+\frac{(q^{(\a-\b)/4}-q^{(\b-\a)/4})q
^{(2n+\a+\b)/4}}{
i(1-q^{n+1+(\a+\b)/2})(1-q^{n+2+(\a+\b)/2})}\right]\,s_n^{(\a,\b)
}(x)} \\
& & -\;\frac{(1-q^{n+\a+1})(1-q^{n+\b+1})q^{n+1+(\a+\b)/2}}{
(1-q^{n+1+(\a+\b)/2})(q^{n+2+(\a+\b)/2};q^{1/2})_3}s_{n-1}^{(\a,\b)}(x).
\nonumber
\eea
We also have the initial conditions
\be
s_0^{(\a,\b)}(x):=1,\qquad s_1^{(\a,\b)}(x):=0.
\ee
When
\be
\a=\overline{\b},\qquad \mbox{Im}\;\a\ne 0,\qquad  \mbox{Re} \; \a
> -1, 
\ee
then the  coefficient of $s_n^{(\a,\b)}(x)$ in (5.19) is real for
$n\ge 0$
and the coefficient of $s_{n-1}^{(\a,\b)}(x)$ is negative for
$n>0$. Thus
the $s_n$'s are orthogonal with respect to a positive measure, say
$d\psi$.
The coefficients in the recurrence relation (5.19) are bounded.
Thus we can
apply Markov's theorem \cite{Sz}, namely
\be
\lim_{n\to\infty}(s_n^{(\a,\b)}(x))^*/s_n^{(\a,\b)}(x)=
\int_{-\infty}^\infty \frac{d\psi(t)}{x-t},\quad\mbox{Im}\;x\ne 0,
\ee
where $(s_n^{(\a,\b)}(x))^*$ is a solution to (5.20) satisfying the
initial
conditions
\be
(s_0^{(\a,\b)}(x))^*:=0,\qquad (s_1^{(\a,\b)}(x))^*:=1.
\ee
It is easy to see that
\be
(s_n^{(\a,\b)}(x))^*=s_{n-1}^{(\a+1,\b+1)}(x).
\ee
Therefore (5.2), (5.18) and (5.23) give
\bea\lefteqn{
\int_{-\infty}^\infty \frac{d\psi(t)}{x-t}=
\frac{(p^{\a+\b+2};p)_2}{(1-p^{\b+1})(1-ip^{\a+3/2}/x)}} \\
& &
\cdot\frac{{}_2\phi_1(p^{\a+2},-ip^{1/2}/x;ip^{\a+5/2}/x;p,p^{\b+
2})}{
{}_2\phi_1(p^{\a+1},-ip^{1/2}x;ip^{\a+3/2}/x;p,p^{\b+1})}.
\nonumber
\eea

Recall that the Coulomb wave function $F_L(\eta,\rho)$ are defined
in terms of 
a confluent hypergeometric function as
\cite{Ab:St}
\be
F_L(\eta,\rho):=2^{L}e^{-\pi\eta/2}\frac{|\G(L+1+i\eta)|}{\G(2L+2)}
\rho^{L+1}e^{-i\rho}{}_1F_1(L+1-i\eta;2L+2;2i\rho)
\ee
As $q\to 1^-$ the polynomials $(1-q)^ns_n^{(\a,\b)}(x/(1-q))$ tend
to the
Wimp polynomials, \cite{Wi}. The right hand side of (5.2) in the
case of
Wimp's polynomials is $F_L(\eta,\rho)$. This suggests defining a
$q$-analog
of $F_L(\eta,\rho)$ by
\be
F_L(\eta,\rho ;q):=(iq^{1/2}\rho; q)_\infty\,{}_2\phi_1
\left(\left.\ba{c}
-q^{L+i\eta+1},q^{L-i\eta+1} \\
q^{2L+2} \ea \right|q,iq^{1/2}\rho\right).
\ee
where $L$ and $\eta$ are real parameters. Observe that the iterate
of the
Heine transformation \cite[(III.3)]{Ga:Ra}
\be
{}_2\phi_1(a,b;c;q,z)=\frac{(abz/c;q)_\infty}{(z;q)_\infty}
{}_2\phi_1(c/a,c/b;c,q,abz/c).
\ee
shows that $F_L(\eta,\rho ;q)$ is real when $\rho$ is real, as in
the case
of  $F_L(\eta,\rho)$.

\newpage
\section{An Expansion Formula.}
\setcounter{equation}{0}
The purpose of this section is to give a direct proof of the 
eigenfunction expansion (6.13).

Set 
\be
{\cal E}_q(x;a,r)=\sum_{m=0}^\infty a_mp_m(x;b,b\sqrt{q},-c,-c
\sqrt{q}).
\ee
In order to compute $a_m$ we need to evaluate the integrals
\be
J_m(a;r):=\Int w(x;b,b\sqrt{q},-c,-c \sqrt{q})
p_m(x;b,b\sqrt{q},-c,-c \sqrt{q}){\cal E}_q(x;a,r)dx,
\ee
when $a=-i$. We shall keep the parameter $a$ in (6.2) free till the

end then we specialize the result by choosing $a=-i$. It is clear
that
\be
J_m(a;r)=\sum_{n=0}^\infty \frac{q^{n^2/4}r^n}{(q;q)_n}
I_{m,n}(a,b,c),
\ee
where
\be
I_{m,n}(a,b,c):=\Int w(x;b,b\sqrt{q},-c,-c \sqrt{q})
p_m(x;b,b\sqrt{q},-c,-c
\sqrt{q})\frac{h(x;aq^{\frac{1-n}2})}{h(x;aq^{\frac{1+n}2})}
dx.
\ee
Formulas (6.3.2) and (6.3.9) in \cite{Ga:Ra} imply
\bea\lefteqn{
\Int w(x;\a,\b,\g,\delta)\frac{h(x;g)
}{h(x;f)}dx }\nonumber\\
& & =\frac{2\pi(\a g,\b g,\g g, \delta g,fg,\a\b\g\delta
f/g;q)_\infty}{
(q,\a\b,\a\g,\a\delta,\a f,\b\g,\b\delta,\b f,\g\delta,\g f,\delta
f,g^2;q)_\infty}
\nonumber\\
& & \quad\cdot{}_8W_7(g^2/q;
g/\a,g/\b,g/\g,g/\delta,g/f;q,\a\b\g\delta f/g).
\nonumber
\eea
Using the $_4\phi_3$ representation of $p_n$ in (6.4) we obtain
\bea\lefteqn{
I_{m,n}(a,b,c)=\frac{2\pi(abq^{(1-n)/2},abq^{1-n/2},-acq^{(1-n)/2},
-acq^{1-n/2};q)_\infty}{(q,b^2\sqrt{q},-bc,-bc\sqrt{q},-bc\sqrt{q},
-bc{q},c^2\sqrt{q},abq^{(1+n)/2};q)_\infty} }\\
& & \quad\cdot\frac{(qa^2,b^2c^2q^{n+1};q)_\infty}{
(abq^{1+n/2},-acq^{(n+1)/2},-acq^{1+n/2},a^2q^{1-n};q)_\infty}
\nonumber \\
& & \cdot \sum_{j=0}^{m} \frac{(q^{-m},b^2c^2q^m,
abq^{(n+1)/2})_j\,q^j}{
(q,b^2c^2q^{n+1},abq^{(1-n)/2})_j}
\nonumber \\
& &\quad\cdot_8W_7(a^2q^{-n};
aq^{-j-(n-1)/2}/b,aq^{-n/2}/b,-aq^{(1-n)/2}/c,
-aq^{-n/2}/c,q^{-n};q,b^2c^2q^{j+n+1}).
\nonumber
\eea
We now apply Watson's transformation formula which expresses a
terminating
very well-poised $_8\phi_7$ as a multiple of a terminating balanced
$_4\phi_3$, \cite[ (III.17)]{Ga:Ra}. Thus
\bea\lefteqn{
_8W_7(a^2q^{-n};aq^{-j-(n-1)/2}/b,aq^{-n/2}/b,-aq^{(1-n)/2}/c,
-aq^{-n/2}/c,q^{-n};q,b^2c^2q^{j+n+1})}\nonumber \\
& & =\frac{(a^2q^{1-n},c^2
q^{1/2})_n}{(-acq^{(1-n)/2},-acq^{(2-n)/2})_n}
\,{}_4\phi_3\left(\left.
\ba{c} q^{-n},-aq^{-n/2}/c,-aq^{(1-n)/2}/c, b^2q^{j+1/2} \\
q^{-n+1/2}/c^2, abq^{j+(1-n)/2}, abq^{1-n/2}\ea \right|\,
q,q\right).\nonumber
\eea
We then apply the Sears transformation (2.4) with invariant
parameters
$q^{-n}, -aq^{-n/2}/c, q^{-n+1/2}/c^2$. After some simplification
we 
obtain
\bea\lefteqn{
_8W_7(a^2q^{-n};aq^{-j-(n-1)/2}/b,aq^{-n/2}/b,-aq^{(1-n)/2}/c,
-aq^{-n/2}/c,q^{-n};q,b^2c^2q^{j+n+1})} \\
& & =\frac{(a^2q^{1-n},c^2 q^{1/2}, -q^{-n}/bc,-q^{-n+1/2}/bc)_n}{
(-acq^{(1-n)/2},-acq^{1-n/2},-abq^{(1-n)/2},-abq^{1-n/2})_n}
\nonumber \\
& & \cdot\frac{(-bcq^{n+1/2}, abq^{(1-n)/2})_j}{
(abq^{(n+1)/2}, -bcq^{1/2})_j}\; (-acb^2q^{(n+1)/2})^j \nonumber \\
& & \cdot  _4\phi_3\left(\left.
\ba{c} q^{-n},q^{-n-j}/b^2c^2, -aq^{-n/2}/c,-q^{-n/2}/ac
\nonumber\\
-q^{-n}/bc, -q^{-n-j+1/2}/bc,q^{-n+1/2}/c^2\ea  \right|
q,q\right).\nonumber
\eea
The substitution of the right-hand side of (6.6) for the $_8W_7$ in
(6.5) and some simplification  lead to
\[
I_{m,n}(a,b,c)=
\kappa(b,c)\frac{(c^2q^{1/2},
-bcq^{1/2},-bcq)_n}{(qb^2c^2)_n}(-a/c)^nq^{-n^2/2}
\]
\dis\dis
\[
\mbox{\hspace{0.4in}}\cdot \sum_{j=0}^m\frac{(q^{-m},b^2c^2q^m,
-bcq^{n+1/2})_j}{(q,b^2c^2q^{n+1},-bcq^{1/2})_j}\,q^j \]
\dis\dis
\[
\mbox{\hspace{0.7in}}\cdot\,
 _4\phi_3\left(\left.
\ba{c} q^{-n},q^{-n-j}/b^2c^2, -aq^{-n/2}/c,-q^{-n/2}/ac \\
-q^{-n}/bc, -q^{-n-j+1/2}/bc,q^{-n+1/2}/c^2\ea  \right| q,q\right),
\]
where 
\be 
\kappa(b,c)=\frac{2\pi (qb^2c^2)_\infty}{(q,b^2q^{1/2},
-bc,-bcq^{1/2}, 
-bcq^{1/2},-bcq,
c^2q^{1/2})_\infty}
\ee
\dis\dis
\[
\mbox{\hspace{0.4in}}=\frac{2\pi (bcq^{1/2},bcq)_\infty}{
(q,b^2q^{1/2}, -bc,-bcq^{1/2},c^2q^{1/2})_\infty}.
\]
Replace the $_4\phi_3$ by its series definition with summation
index $k$.
Then interchange the $j$ and $k$ sums to obtain
\[
I_{m,n}(a,b,c)=
\kappa(b,c)\frac{(c^2q^{1/2},
-bcq^{1/2},-bcq)_n}{(qb^2c^2)_n}(-a/c)^nq^{-n^2/2}
\]
\dis\dis
\[
\mbox{\hspace{0.4in}}\cdot \sum_{k=0}^n\frac{(q^{-n},-aq^{-n/2}/c,
-q^{-n/2}/ac,q^{-n}/b^2c^2)_k}{(q,-q^{-n}/bc, q^{-n+1/2}/c^2,
-q^{-n+1/2}/bc)_k}\,q^k
\]
\dis\dis
\[
\mbox{\hspace{0.7in}}\cdot{}_3\phi_2\left(\left. \ba{c}
q^{-m}, b^2c^2q^m,-bcq^{n-k+1/2} \\
b^2c^2q^{n+1-k},-bcq^{1/2}\ea  \right|\,q,q\right).
\]
The ${}_3\phi_2$ can now be summed by the $q$-analog of the
Pfaff-Saalsch\"
utz theorem \cite[(II.12)]{Ga:Ra}. It's sum is
\[
(q^{k-n},q^{-m+1/2}/bc)_m/(-bcq^{1/2},q^{k-m-n}/b^2c^2)_m,
\]
which vanishes for all $k,\; 0\le k \le n$ if $m>n$. Thus we get
\be
I_{m,n}(a,b,c)=
\frac{(q^{-n},-q^{-m+1/2}/bc)_m}{(-bcq^{1/2},q^{-m-n}/b^2c^2)_m}
\kappa(b,c)\frac{(c^2q^{1/2},
-bcq^{1/2},-bcq)_n}{(qb^2c^2)_n}(-a/c)^nq^{-n^2/2}
\ee
$$\qquad
\mbox{\hspace{0.4in}}\cdot{}_4\phi_3\left(\left. \ba{c}
q^{m-n},q^{-m-n}/b^2c^2,-aq^{-n/2}/c,
-q^{-n/2}/ac\\
-q^{-n}/bc,-q^{-n+1/2}/bc, q^{-n+1/2}/c^2\ea \right|\,q,q\right).
$$
The relationship (6.3) and the observation $I_{m,n}(a,b,c)=0$ if
$n<m$ show
that $J_m(a;r)$ is given by
\[
\frac{q^{m^2/4}r^m}{(q)_m}\sum_{n=0}^\infty
q^{n^2/4}\frac{(rq^{m/2})^n}{
(q^{m+1})_n}\,I_{n,n+m}(a,b,c).
\]
Applying (6.8) we find after some simplification
\be
J_m(a,r)=\kappa(b,c)q^{m^2/4}(-abr)^m\frac{(c^2q^{1/2}, -bcq^{1/2},
-bcq)_m}{(qb^2c^2)_{2m}}
\ee
\dis\dis
\[
\mbox{\hspace{0.4in}}\cdot \sum_{n=0}^\infty
\frac{(c^2q^{m+1/2}, -bcq^{m+1/2},
-bcq^{m+1})_n}{(q,b^2c^2q^{2m+1})_{n}}
\left(-\frac{ar}c\,q^{-m/2}\right)^n
\]\dis\dis
\[
\mbox{\hspace{0.4in}}\cdot  q^{-n^2/4}{}_4\phi_3\left(\left. \ba{c}
q^{-n},q^{-2m-n}/b^2c^2,-aq^{-(n+m)/2}/c, -q^{-(n+m)/2}/ac \\
-q^{-n-m}/bc, -q^{-n-m+1/2}/bc,q^{-n-m+1/2}/c^2 \ea
\right|\,q,q\right).
\]
The ${}_4\phi_3$ on the right hand side of (6.9) has a quadratic
transformation. By (3.10.13) of \cite{Ga:Ra}, the aforementioned
${}_4\phi_3$ 
 is equal to
\[{}_4\phi_3\left(\left. \ba{c}
q^{-n/2},q^{-m-n/2}/bc,-aq^{-(n+m)/2}/c, -q^{-(n+m)/2}/ac \\
-q^{-n-m}/bc, -q^{(-n-m+1/2)/2}/c, q^{(-n-m+1/2)/2}/c \ea
\right|\,q^{1/2},q^{1/2}\right)
\]
\dis\dis
\[
\mbox{\hspace{0.4in}}=\frac{(q^{-n/2},q^{-m-n/2}/bc,-aq^{-(n+m)/2
}/c,
-q^{-(n+m)/2}/ac;q^{1/2})_n}{(q^{1/2},-q^{-n-m}/bc,
q^{(-n-m+1/2)/2}/c,-q^{(-n-m+1/2)/2}/c;q^{1/2})_n}\,q^{n/2}
\]
\dis\dis
\[
\mbox{\hspace{0.4in}}\cdot {}_4\phi_3\left(\left. \ba{c}
q^{-n/2}, -bcq^{m+(n+1)/2},cq^{(m+1/2)/2}, -cq^{(m+1/2)/2} \\
bcq^{m+1/2}, -cq^{(m+1)/2}/a,-acq^{(m+1)/2} \ea
\right|\,q^{1/2},q^{1/2}\right),
\]
by reversing the sum in the first ${}_4\phi_3$. After some
straightforward manipulations we establish
\be
J_m(a,r)=\kappa(b,c)\frac{(c^2q^{1/2};
q)_m(-abr)^m}{(bcq^{1/2},bcq;q)_{m}}q^{m^2/4}
\ee
\dis\dis
\[
\mbox{\hspace{0.4in}}\cdot \sum_{n=0}^\infty
\frac{(-cq^{(m+1)/2}/a, -acq^{(m+1)/2};q^{1/2})_n}{
(q^{1/2}, -q^{1/2};q^{1/2})_n}
\left(-\frac{ar}c\,q^{-m/2-1/4}\right)^n
\]\dis\dis
\[
\mbox{\hspace{0.4in}}\cdot{}_4\phi_3\left(\left. \ba{c}
q^{-n/2}, -bcq^{m+(n+1)/2},cq^{(m+1/2)/2}, -cq^{(m+1/2)/2}/ac \\
bcq^{m+1/2}, -cq^{(m+1)/2}/a,-acq^{(m+1)/2} \ea
\right|\,q^{1/2},q^{1/2}\right).
\]
Finally we apply Sears transformation (2.4) to the ${}_4\phi_3$ in
(6.10) with invariant parameters $q^{-n}$, $cq^{(m+1/2)/2}$ and
$bcq^{m+1/2}$. This enables us to cast (6.10) in the form
\be
J_m(a,r)=\kappa(b,c)\frac{(c^2q^{1/2};
q)_m(-abr)^m}{(bcq^{1/2},bcq;q)_{m}}q^{m^2/4}
\ee
\dis\dis
\[
\mbox{\hspace{0.4in}}\cdot \sum_{n=0}^\infty
\frac{(-aq^{1/4}, -q^{1/4}/a;q^{1/2})_n}{
(q^{1/2}, -q^{1/2};q^{1/2})_n}
(ar)^n
\]\dis\dis
\[
\mbox{\hspace{0.4in}}\cdot{}_4\phi_3\left(\left. \ba{c}
q^{-n/2},-q^{-n/2},cq^{(m+1/2)/2}, -bq^{(m+1/2)/2}    \\
bcq^{m+1/2}, -aq^{(-n+1/2)/2}, -q^{(-n+1/2)/2}/a   \ea
\right|\,q^{1/2},q^{1/2}\right).
\]
It is evident from (6.11) that $J_m(a;r)$ is a double series. When
$a^2=-1$ we have been able to reduce the right-hand  side of (6.11)
to
a single series. To see this, replace the ${}_4\phi_3$ in (6.11) by

its defining series then interchange the sums. The result is
\[
J_m(-i;r)=\kappa(b,c)\frac{(c^2q^{1/2};
q)_m}{(bcq^{1/2},bcq;q)_{m}}\frac{(irq^{1/2};q)_\infty}{
(-ir;q)_\infty}(ibr)^mq^{m^2/4}
\]\dis\dis
\[
\mbox{\hspace{0.4in}}\cdot{}_2\phi_1\left(\left. \ba{c}
cq^{m/2+1/4}, -bq^{m/2+1/4}\\
bcq^{m+1/2} \ea \right|\,q^{1/2},ir\right).
\]
We have tried to express  $J_m(a;r)$  of (6.11) as a single sum for
general $a$ but this does not seem to be possible except when
$a=\pm i$.

Now the orthogonality relation (1.10) gives for the case
$a=-i$
\be
\kappa(b,c)a_m=J_m(-i;r)\frac{(1-b^2c^2q^{2m})(b^2c^2,b^2q^{1/2},
-bc;q)_m}
{(1-b^2c^2)(q,c^2q^{1/2},-bcq;q)_m}\,b^{-2m}.
\ee
Thus we established the expansion formula
\be
{\cal E}_q(x;-i,r)=\sum_{m=0}^\infty a_m
p_m(x;b,bq^{1/2},-c,-cq^{1/2}),
\ee
with the $a_m$'s given by
\be
a_m=\frac{(b^2c^2,b^2q^{1/2};q)_m\; (irq^{1/2}; q)_\infty}
{(q,bcq^{1/2},bc;q)_m \; (-ir; q)_\infty}
(ir/b)^mq^{m^2/4} \mbox{\hspace{0.06in}}{}_2\phi_1\left(\left.
\ba{c}
cq^{m/2+1/4}, -bq^{m/2+1/4}\\
bcq^{m+1/2} \ea \right|\,q^{1/2},ir\right).
\ee

\newpage
\section{A Formal Approach}
\setcounter{equation}{0}
We now  formalize the procedure followed  in Section 3
 and used earlier in \cite{Is:Zh}.  Let  $S \subset {\cal C}^k$ and
assume that for every $A \in S$, the sequence of polynomials
$\{p_n(x; A)\}_0^\infty$ are orthogonal with respect to  a measure
with a
nontrivial  absolutely continuous component.  By $A + 1$ we mean  
$(1 + a_1, \dots, 1 +  a_k)$ if 
$A = (a_1, \dots , a_k)$. We will assume that $A + 1 \in S$
whenever $A \in S$.
Let the orthogonality relation of the $p_n$'s be
\be
\int_{-\infty}^{\infty} p_n(x; A) p_m(x; A) d\mu(x; A) = h_n(A)
\;\delta_{m,n}.
\ee
Assume that $\cal D$ is an operator  defined on polynomials by
linearity and
by its action on the basis $\{p_n(x; A)\}$ via  
\be
{\cal D} p_n(x; A) = \xi_n(A) \; p_{n-1}(x; A + 1).
\ee
Furthermore assume that  the support of $\mu'(x;A) = \frac{d\mu(x;
A)}{dx}$ is the same
 for all $A \in S$ and that we know the connection coefficients in
\be
p_n(x; A) = \sum_{j = 0}^n c_{n,j}\; p_j(x; A+1).
\ee
Therefore 
\be
c_{n,j} =  \frac{1}{h_j(A + 1)} \int_{-\infty}^{\infty}  p_n(x;
A)\, p_j(x; A+1)
\; d\mu(x; A+1).
\ee
The formula dual to (7.3) is
\be
 p_n(x; A+1) \; \mu'(x; A+1) = \sum_{m = n}^{\infty} 
\frac{h_{n}(A+1)}{h_m(A)}
\; c_{m,n}p_m(x; A)\; \mu'(x; A),
\ee
holding on the interior of the support of $\mu'(x; A)$.

We now wish to describe the spectrum of a formal inverse to $\cal
D$. 
 Note that 
we can define $\cal D$ densely on $L^2(d\mu(x; A))$ by (7.2) 
provided that the polynomials $p_n(x; A)$ are dense in 
$L^2(d\mu(x; A))$ for all $A \in S$.  This suggests that we define
${\cal D}^{-1}$
via
\be
{\cal D}^{-1} \sum_{n = 0}^\infty a_n\, p_n(x; A + 1) := 
\sum_{n = 0}^\infty a_n\, p_{n+1}(x; A)/\xi_{n+1}(A).
\ee
This motivates the definition
\be
(T_Af)(x) := \int_{-\infty}^\infty f(t)\left[\sum_{n = 0}^{\infty}
\frac{p_{n+1}(x; A)\, p_{n}(t; A+1)} {\xi_{n+1}(A) \,
h_n(A+1)}\right]
 d\mu(t; A+1),
\ee
if $f \in L^2(d\mu(x; A+1))$.

The next step is to consider the eigenvalue problem
\be
T_A\, g = \l\, g, \quad  g(x) \approx \sum_{n = 0}^\infty a_n(\l;A)
\, 
p_{n}(x; A).
\ee
In order for (7.8) to hold it is necessary that $g$ lies in the
range of $T_A$, 
hence (7.6) shows that $a_0(\l;A) = 0$.  Therefore (7.7) and (7.8)
yield
\bea
\sum_{n = 1}^\infty a_n(\l;A) \, 
p_{n}(x; A) = \sum_{n = 1}^\infty \frac{p_{n+1}(x;
A)}{\xi_{n+1}(A)} \, 
  \sum_{k = n}^\infty  \frac{a_k(\l;A)}{  
h_{n}(A+1)} \; \int_{-\infty}^{\infty} p_{k}(t; A) p_{n}(t; A+1)
d\mu(t; A+1). \nonumber
\eea
Thus we have established
\be
\l \; \sum_{n = 1}^\infty a_n(\l;A) \, 
p_{n}(x; A) = \sum_{n = 1}^\infty \frac{p_{n+1}(x;
A)}{\xi_{n+1}(A)} \, 
\sum_{k = n}^{\infty} c_{k,n}\; a_k(\l,A).
\ee
Now (7.9)  implies the recurrence relation
\be
\l \, \xi_n(A) \, a_n(\l; A) = \sum_{k = n-1}^{\infty} c_{k,n-1}\,
a_k(\l, A).
\ee

Observe that (7.10) transformed the eigenvalue problem (7.8) to the
discrete
eigenvalue problem (7.10). When $c_{n,k}  = 0$ for $k < n - r$ for
a fixed $r$
then (7.10) is the eigenvalue equation of a matrix with at most 
$r+1$ nonzero entries in each row. The cases analyzed in
\cite{Is:Zh} 
and in this paper  are the cases when
\be
c_{n,k}  = 0\quad for \; k < n - 2, \; n = 2, 3, \cdots.
\ee
Note that $c_{n,n} \ne 0$. When (7.11) holds then (7.10) reduces to

\be
\l \, \xi_n(A) \, a_n(\l; A) =c_{n-1,n-1} \,a_{n-1}(\l; A)\; +
c_{n,n-1}
a_n(\l; A) \; + c_{n+1,n-1} a_{n+1}(\l; A).
\ee
For example in the case when  the $p_n$'s are the ultraspherical
polynomials
$C_n^\nu(x)$ we have \cite[, \S 144]{Rai},
\be
\frac{d}{dx} C_n^\nu(x) = 2\nu C_{n-1}^{\nu + 1}(x), \quad 2(n+\nu)
C_n^\nu(x)
= 2\nu [C_{n}^{\nu + 1}(x) - C_{n-2}^{\nu + 1}(x)].
\ee
Therefore 
\bea
\xi_n(\nu) = 2\nu,\; c_{n,n} = \quad -c_{n, n-2} = \nu/(\nu + n),
\quad 
c_{n, n-1} = 0, \nonumber
\eea
and (7.12) becomes
\be
2\,\l \, a_n(\l;\nu) = \frac{a_{n-1}(\l;\nu)}{(\nu +n-1)} -\; 
\frac{a_{n+1}(\l;\nu)}{(\nu +n+ 1)},
\ee
which is (2.11) in \cite{Is:Zh}.

The procedure just outlined is very formal but can be justified 
if both $p_n(x; A)$ and $ p_n(x; A+1)$ are dense in $L^2(d\mu(x;
A))
\cap L^2(d\mu(x; A+1))$.

Another approach to the same problem is to think of $T_A$ as a
right inverse to
 ${\cal D}$. In other words ${\cal D}T_A$ is the restriction of the
identity   
operator to the range of ${\cal D}$. Thus (7.8) is equivalent to
\be
g(x; \l) = \l {\cal D} g(x; \l).
\ee
Now the use of the orthogonal expansion of $g$ and formulas (7.2) 
and (7.15) implies  that $\l\, \xi_{n+1}(A) \; a_{n+1}(\l, A)$ is
the projection 
of $\sum_{k=1}^{\infty} a_k(\l;A)\; p_k(x; A)$ on the space spanned
by
$p_n(x; A+1)$. Therefore (7.3) implies (7.12).

In \cite{Is:Zh} it was observed that the eigenfunction expansion 
$\sum_{k=1}^{\infty} a_k(\l;A)\; p_k(x; A)$ can be extended to
values of $\l$
off the discrete spectrum of the operator under consideration. 
This is also the case with the expansion formula (6.13). We now
attempt to find 
such an expansion in general.

We seek functions $\{F_n(\l;A)\}$ such that the function$E_A(x;
\l)$,
\be
E_A(x; \l) := \sum_{k=0}^{\infty} F_k(\l;A)\; p_k(x; A),
\ee 
satisfies
\be
\l\, {\cal D} E_A(x; \l) = E_A(x; \l).
\ee
Observe that (7.16) reminds us of (7.15), hence under the
assumption (7.11),
(7.3) and (7.17) imply that
the $F_n$'s must satisfy a recursion relation similar to (7.12),
that is
\be
\l\,\xi_n(A) F_n(\l;A) = c_{n-1,n-1} \,F_{n-1}(\l;A) +
c_{n,n-1}\,F_n(\l;A)
+c_{n+1,n-1}\,F_{n+1}(\l;A).
\ee
In order to maintain a parallel course with the results in
\cite{Is:Zh} 
and with the notation of Section 5 we will renormalize the $a_n$'s
and $F_n$'s
in order to put (7.12) in monic form and change (7.18) to a
recursion with the 
coefficient of $y_{n-1}$ equal to unity.  Keeping in mind that
$a_0(x;\l) = 0$
and that $a_1(x; \l)$ is a multiplicative constant which we may
take to be
unity, we set
\be
a_n(\l;A) = \prod_{j = 0}^{n-2} \frac {\xi_{j+1}(A)} {c_{j+2,j}}\;
 u^{n-1}\; b_{n-1}(\l u; A)\; n > 1, \quad  a_1(x; \l) = b_0(x; \l)
= 1.
\ee
and
\be
F_n(\l;A) = \prod_{j = 0}^{n-1} \frac{c_{j,j}} {\xi_{j+1}(A)}\;
 u^{n-1}\; G_{n-1}(\l u; A), \; n > 0\quad F_0(\l;A) = G_{-1}(\l u;
A).
\ee
Here $u$ is a free normalization factor at our disposal and  may
depend 
on $A$.  Thus 
\be
\l b_n(\l; A) = \; b_{n+1}(\l; A) + \; B_n(A)  b_n(\l; A) + \;
C_n(A) b_{n-1}(\l; A),
\ee
and
\be
\l G_n(\l; A) = \; C_{n+1}(A)G_{n+1}(\l; A) + \; B_n(A)  G_n(\l; A)
+ \; 
 G_{n-1}(\l; A),
\ee
hold with
\be
B_n(A) := \frac{u\,c_{n+1,n}}{\xi_{n+1}(A)},\quad C_n(A) =
 \frac{u^2\, c_{n,n}\; c_{n+1,n-1}}{\xi_{n}(A)\, \xi_{n+1}(A)}.
\ee

We are seeking a solution to (7.18) that makes (7.16) converge on
sets of $\l$'s
containing the spectrum. Since the $p_n$'s are given we need to
choose the
$F_n$'s to be as small as possible, that is choose $F_n$ to be the
minimal 
solution, if  it exists. Recall that a solution  $w_n$ of (7.18) is
minimal if
 $w_n = o(v_n)$ where $v_n$ is any other linearly independent
solution of 
the same recurrence relation, \cite{Jo:Th}. It is clear that a
minimal solution of (7.18) exists if and only if (7.22) has a
minimal solution. The minimal solution
may change form in different regions of the parameter or variable
space, 
\cite{Gu:Is}.
 When $B_n(A) \to 0$ and
 $C_n(A) \to 0$ as $n \to \infty$  then
(7.22) has a minimal solution, see Theorem 4.55 in \cite{Jo:Th}.

In many cases we encounter a fortutious  situation where we can
choose $u$ 
in (7.23) such that 
\be
B_n(A) = B_0(A+n),\quad C_n(A) = C_0(n+A).
\ee
Therefore
\be
G_n(\l;A) = G_0(\l;A+n) = G_{-1}(\l;A+n+1).
\ee
When (7.24) holds Theorem 4.5 of \cite{Is:Zh} comes in handy. This
latter
 theorem is
\begin{thm}
Let $f(x;\nu + A)$ be a multi-parameter family of functions
satisfying
\be
C_{\nu + A} f(x;\nu + A+1)=(A_{\nu + A} x+B_{\nu + A})f(x;\nu + A)
\pm\;f(x;\nu + A-1),
\ee
and let $\{f_{n,\nu + A}(x)\}$ be a sequence of polynomials defined
by
\be
f_{0,\nu + A}(x)=1,\quad f_{1,\nu + A}(x)=A_{\nu + A} x+B_{\nu +
A},
\ee\dis\dis
\be
f_{n+1,\nu + A}(x)=(A_{n+\nu + A} x+B_{n+\nu + A})f_{n,\nu +
A}(x)\;\pm\,C_{n+\nu + A-1}f_{n-1,\nu + A}(x).
\ee
Then
\be
C_{\nu + A} C_{\nu + A+1}\cdots C_{\nu + A+n-1} f(x;\nu +
A+n)=f_{n,\nu + A}(x) f(x,\nu + A)\pm
f_{n-1,\nu + A+1}(x) f(x;\nu + A-1).
\ee
\end{thm}
Theorem 7.1  establishes
\bea \lefteqn{
C_{-1}(A)\, C_{-1}(A+1)\cdots C_{-1}(A+n-1)G_0(\l; A+n)}\\
& & = b_n(\l;A) G_0(\l;A)
+ b_{n-1}(\l;A+1) G_{-1}(\l;A),\nonumber
\eea
where we used $G_{-1}(\l;A) = G_{0}(\l;A-1)$, see (7.25).

Now assume  in addition to (7.26) that  $B_n(A) \to 0$ and $C_n(A)
\to 0$,
 hence the minimal solution to (7.22) exists. Let $\{G_n(\l;A)\}$ 
be the minimal solution to (7.22).
 According to Pincherle's theorem, \cite{Jo:Th}, the continued
$J$-fraction 
associated with (7.21) converges to a constant multiple of
 $G_{0}(\l;A)/G_{-1}(\l;A)$. Therefore the eigenvalues of the
infinite 
tridiagonal matrix associated with (7.21) are the zeros of
$G_{-1}(\l;A)$ 
which are not zeros of $G_{0}(\l;A)$. If  $G_{-1}(\l;A) = 0$ then
(7.30) indicates that the series in (7.16)  is a multiple of the
series (7.8) 
and the multiplier does not depend on $x$ but may depend on $\l$. 
This 
explains the relationship between the expansion representing the 
eigenfunction $g$ in (7.8) and the expansion in (7.16) which is
expected
to be valid for a range of $\l$ wider than the spectrum of ${\cal
D}$.

Finally we apply the preceding outline to the case of continuous
 $q$-Jacobi polynomials and give another proof of (6.13)-(6.14). 
In the case under consideration
\be
\xi_n(\a, \b) := \frac {2q^{-n + (\a + 5/2)/2} \,(1 - q^{n+\a + \b
+ 1})}
{(1-q)\, (-q^{(\a+\b+1)/2}; q^{1/2})_2},
\ee
\be
c_{n,n} := \frac{q^{-n/2}(1 - q^{\a + \b + n + 1})
(1 - q^{\a + \b + n + 2})}{(- q^{(\a + \b  + 1)/2}; q^{1/2})_2(1 -
q^{n + (\a + \b + 1)/2})(1 - q^{n + (\a + \b + 2)/2})},
\ee
\be
c_{n,n-1} :=\frac{q^{(\a + \b + 2 -n)/2}(1 - q^{\a + \b + n + 1})
(1 + q^{n + (\a + \b  + 1)/2}) (1 - q^{(\a - \b)/2})}
{(- q^{(\a + \b  + 1)/2}; q^{1/2})_2(1 - q^{n + (\a + \b)/2})(1 -
q^{n + (\a + \b + 2)/2})}, 
\ee
and
\be
c_{n,n-2} :=
- \frac{q^{(3\a + \b + 4 -n)/2}(1 - q^{\a + n})(1 - q^{\b + n})}
{(- q^{(\a + \b  + 1)/2}; q^{1/2})_2(1 - q^{n + (\a + \b)/2})(1 -
q^{n + (\a + \b + 1)/2})}.
\ee
With the choice
\be
u = \frac{2q^{1/2}}{1-q}
\ee
we find
\be
B_n(\a,\b) = B_0(\a+n,\b+n) = - \frac{(1-q^{\frac{\b-\a}2})
(1+q^{\frac{\a+\b+3}2+n})}{(1-q^{\frac{\a+\b+2}2+n})
(1-q^{\frac{\a+\b+4}2+n})}q^{(n + \a + 3/2)/2}
\ee
and
\be
 C_n(\a,\b) = C_0(\a+n,\b+n) = - \frac{(1-q^{\a+1+n})
(1-q^{\b+1+n})q^{n+\frac{\a+\b}{2}+1}}
{(1-q^{\frac{\a+\b+1}2+n})(1-q^{\frac{\a+\b+2}2+n})^2
(1-q^{\frac{\a+\b+3}2+n})}
\ee
Therefore (5.13) yields
\be
G_n(\l; \a, \b) = G_0(\l;  \a +n, \b+n) =
(-\l)^{-(\a+\b)/2}X_n^{(\a,\b)}(\l).
\ee
When  $A = (\a, \b)$ and $p_n(x; A)$ are the continuous 
$q$-Jacobi polynomials  $P_n^{(\a,\b)}(x|q)$ we will denote $E_A(x;
\l)$ by 
$E_{\a,\b}(x; \l)$.
Now  with ${\cal D} =  {\cal D}_q$ formula (7.16) becomes
\be
E_{\a, \b}(x; \l) = \Summ q^{n(n - 2\a)/4}\frac{(q^{\a + \b + 1};
q)_n}
{(q^{(\a + \b + 1)/2};
q^{1/2})_{2n}}\left(\frac{2q^{1/2}}{1-q}\right)^{n-1}
G_{n-1}\left(\frac{2\l q^{1/2}}{1-q}; \a, \b\right)
P_n^{(\a,\b)}(x|q).
\ee
We now find another solution to (7.17) and prove a uniqueness 
theorem for solutions of (7.17). We then  equate $E_{\a, \b}(x;
\l)$
and the second solutions to (7.17) and establish (6.13).
\begin{lem}
Assume that $f(x)$ is an entire function of the complex variable
$x$. If
\be
 {\cal D}_q\,f(x)=\frac{iyq^{1/4}}{1-q}\,f(x),
\ee
then $f(x)$ is unique up to a multiplicative function of $y$.
\end{lem}
Lemma 7.2 is essentially Lemma 3.5 in \cite{Is:Zh}. 

A calculation gives
\be
{\cal D}_q\,{\cal E}_q(x;a,b) =\frac{-2abq^{1/4}}{1-q}{\cal
E}_q(x;a,b).
\ee
Therefore Lemma 7.2 implies
\begin{thm}
If $f$ satisfies the assumptions in Lemma 7.2 then 
\be
f(x) = w(y)\; {\cal E}_q(x; -i, y/2).
\ee
\end{thm}
\begin{thm}
The function $[-2\l q^{1/2}/(1-q)]^{(\a+\b)/2} \, E_{\a,\b}(x; \l)$
does
 not depend on $\a$ or $\b$.
\end{thm}
{\bf Proof.} In general we have
\bea\lefteqn{
 E_{A}(x; \l) = \Summ F_n(\l;A) \, p_n(x; A)  } \\
& & = \Summ F_n(\l;A) [c_{n,n} p_n(x; A+1) + c_{n,n-1} p_{n-1}(x;
A+1)
+ c_{n,n-2} p_{n-2}(x; A+1)]   \nonumber \\
& & =\Summ p_n(x; A+1) [c_{n,n} F_n(\l; A) + c_{n+1,n} F_{n+1}(\l;
A)
+ c_{n+2,n} F_{n+2}(\l; A)]  \nonumber \\
& &  =\Summ p_n(x; A+1) \l \xi_{n+1}(A) F_{n+1}(\l; A) \nonumber \\
& &= \l u \Summ p_n(x; A+1) \, u^{n-1} \, \prod_{j = 0}^n
\frac{c_{j,j}}{\xi_{j+1}(A)}
\; \xi_{n+1}(A)\, G_n(\l u; A)  \nonumber  \\
& & = \l   \Summ c_{n,n} \, p_n(x; A+1) \, u^{n} \, \prod_{j = 0}^n
\frac{c_{j,j}}{\xi_{j+1}(A)} G_{n-1}(\l u; A+1).  \nonumber 
\eea
In the case of continuous $q$-Jacobi polynomials the last equation
gives
\bea\lefteqn{
E_{\a,\b}(x; \l) = \l \Summ \frac{q^{-n/2}
(1 - q^{\a + \b + n + 1})
(1 - q^{\a + \b + n + 2})}{(- q^{(\a + \b  + 1)/2}; q^{1/2})_2(1 -
q^{n + (\a + \b + 1)/2})(1 - q^{n + (\a + \b + 2)/2})}} \\
& & \cdot q^{n(n - 2\a)/4}\frac{(q^{\a + \b + 1}; q)_n}
{(q^{(\a + \b + 1)/2};
q^{1/2})_{2n}}\left(\frac{2q^{1/2}}{1-q}\right)^{n}
G_{n-1}\left(\frac{2\l q^{1/2}}{1-q}; \a+1, \b+1 \right)
P_n^{(\a,\b)}(x|q).
\nonumber
\eea
Therefore
\bea
E_{\a,\b}(x; \l) = \frac{[-2\l q^{1/2}/(1-q)](q^{\a+ \b +1};q)_2}
{(q^{(\a+ \b +1)/2}, -q^{(\a+ \b +1)/2}; q^{1/2})_2}\; E_{\a + 1,\b
+ 1}(x; \l).
\nonumber
\eea
The above functional equation can be put in the form 
\bea
[-2\l q^{1/2}/(1-q)]^{(\a+\b)/2}  E_{\a,\b}(x; \l) = 
[-2\l q^{1/2}/(1-q)]^{(\a+\b+4)/2}  E_{\a+2,\b+2}(x; \l), \nonumber
\eea
and we have
\bea \lefteqn{
[-2\l q^{1/2}/(1-q)]^{(\a+\b)/2}  E_{\a,\b}(x; \l)}\\
 & & = \lim_{m \to \infty}
[-2\l q^{1/2}/(1-q)]^{2m+(\a+\b)/2}  E_{2m + \a, 2m + \b}(x; \l).
\nonumber 
\eea
Now substitute the right-hand sides of (7.38) and (7.39) for
$G_{n}$ 
and $E_{\a,\b}$ in the right-hand side of (7.45) to get
\bea
[-2\l q^{1/2}/(1-q)]^{2m+(\a+\b)/2}  E_{2m + \a, 2m + \b}(x; \l)
\approx \l \Summ q^{n^2/4}\, (-\l)^{-n} q^{-n\a/2}P_n^{(\a,
\b)}(x|q).
\eea
But (1.24) and (1.31) imply
\bea
P_n^{(\a, \b)}(x|q) \approx \frac{q^{n\a/2}}{(q;q)_n} p_n(x;
0,0,0,0|q) = 
\frac{q^{n\a/2}}{(q;q)_n} H_n(x|q), \nonumber
\eea
where $\{H_n(x|q)\}_0^\infty$ are the continuous $q$-Hermite
polynomials.
Thus the limit on the right-hand side of (7.45) exists and we have
established
\be
[-2\l q^{1/2}/(1-q)]^{(\a+\b)/2}  E_{\a,\b}(x; \l) = 
\Summ \frac{q^{n^2/4}\, (-\l)^{-n}}{(q;q)_n}\; H_n(x|q).
\ee
This proves Theorem 7.4.

\begin{cor}
We have
\be
[-2\l q^{1/2}/(1-q)]^{(\a+\b)/2}  E_{\a,\b}(x; \l) =
(\l^{-2};q^2)_\infty
{\cal E}_q(x; -i, i/\l).
\ee
\end{cor}

Corollary 7.5 follows from \cite{Is:Zh} where Ismail and Zhang
proved 
that the right-hand sides of (7.47) and (7.48) are equal.

\begin{thm}
The expansion  of ${\cal E}_q(x; -i, r)$ in a continuous $q$-Jacobi
series
 is given by (6.13) where $b = q^{(2\a+1)/4}$ and $c = 
q^{(2\b+1)/4}$.
\end{thm}

{\bf Proof}. From Theorem 7.3 and Corollary 7.4 we see that the
right-hand 
side of (6.13) is $w(r) {\cal E}_q(x; -i, r)$. Furthermore $w(r)$
does not depend on $\a$ or $\b$ since neither ${\cal E}_q(x; -i,
r)$ nor the right-hand side of
 (6.13) depend on $\a$ or $\b$. Now $w$ can be found by letting
$\a$ and $\b$
tend to $\infty$ then use (7.48).

\newpage
\section{Remarks.}
\setcounter{equation}{0}
In 1940 Schwartz  published  an interesting paper \cite{Sc}
containing the following result.
\begin{thm}
Let $\{p_{n,\nu}(x)\}$ be a family of monic polynomials generated
by
\be
p_{0,\nu}(x)= 1,\; p_{1,\nu}(x) = x + B_{\nu},
\ee
\be
p_{n+1,\nu}(x) = (x + B_{n + \nu})\;p_{n,\nu}(x) + C_{n+\nu} \;
p_{n-1,\nu}(x).
\ee
If  both
\be
\Summ |B_{n+\nu} - a| \;< \; \infty \quad and  \quad \Summ
|C_{n+\nu}|
\;< \; \infty
\ee
hold, then $x^n\; p_{n,\nu}(a + 1/x)$ converges on compact subsets
of the complex plane to an entire function.
\end{thm}
It is clear from (8.1) and (8.2) that $(p_{n,\nu}(x))^* =
p_{n-1,\nu+1}(x)$,
 hence the continued $J$-fraction associated with (8.1) and (8.2)
converges to a 
meromorphic function of $1/x$ and the convergence is uniform on
compact subsets
of the complex plane which neither contain the origin nor contain
poles of the
limiting function. Schwartz illustrated his theory by applying it
to the
Lommel polynomials and he mentioned their orthogonality relation.

Somehow Schwartz's  interesting paper \cite{Sc} was not noticed 
and neither his 
results were   quoted nor his paper was cited in the standard
modern 
references on orthogonal polynomials  \cite{Ch}, \cite{Fr},
\cite{Sz} 
and continued 
fractions, \cite{Wa}, \cite {Jo:Th}.  Many of Schwartz's results
were
later rediscovered by others.  Dickinson, Pollack and Wannier
\cite{Di:Po} rediscovered the special case
$B_{n+\nu} = 0$ of Schwartz's theorem. It is worth noting that if
$B_{n+\nu} = 0$ then  a theorem of Van Vleck, Theorem 4.55 in
\cite{Jo:Th}, 
states that $C_{n+\nu} \to  0$ 
suffices to establish the uniform convergence of the continued
$J$-fraction
 associated with (8.1) and (8.2)  to a meromorphic function. The
convergence being 
 uniform on compact subsets
of the complex plane which neither contain the origin nor contain
poles of the
limiting function.

In \cite{Gu:Is} it was proved that
\be
X_n^{(5)}(x) := \left(-\frac{D}{xABC}\right)^n \frac{(Dq^{2n},
Dq^{2n-1},
-q^n/x)_\infty}{(Aq^n, Bq^n, Cq^n, Dq^n/A, Dq^n/B, Dq^n/C)_\infty}
\ee
\bea
\qquad \qquad .\mbox{}_3\phi_2\left(\left.\ba{c} Aq^{n}, Bq^{n},
Cq^{n}\\
\quad Dq^{2n},\quad -q^{n}/x, \ea \right|q,-\frac{D}{xABC}\right),
\nonumber
\eea
satisfies the three term recurrence relation
\be
Z_{n+1}(x) = (x - a_n)\;Z_n(x) - \;b_n\; Z_{n-1}(x),
\ee
with 
\be
a_n := -\frac{D}{ABC} -q^{n-1}\frac{(1- Dq^n/A)(1- Dq^n/B)(1-
Dq^n/C)}
{(1 - Dq^{2n-1})(1 - Dq^{2n-2})}
\ee
\bea
\qquad \qquad
+\frac{D}{ABC}\; 
\frac{(1-Aq^{n-1})(1-Bq^{n-1})(1-Cq^{n-1})}{(1 - Dq^{2n-1})(1 -
Dq^{2n})},
\nonumber
\eea
and
\be
b_n := -\frac{D}{ABC}q^{n-2}(1-Aq^{n-1})(1-Bq^{n-1})(1-Cq^{n-1})
\ee
\bea
\qquad \qquad
.\frac{(1-Dq^{n-1}/A)(1-Dq^{n-1}/B)(1-Dq^{n-1}/C)}
{(1 - Dq^{2n-1})(1 - Dq^{2n-2})^2(1 - Dq^{2n-3})}. \nonumber
\eea
We next identify (3.12) as a limiting case  $A \to \infty$ of (8.5)
with $B = - C$. When $A \to \infty$, it is easy to see that 
\be
a_n \to -\frac{q^{n-1} (1 + Dq^{2n-1})\, (1 - D/B^2)}
{(1 - Dq^{2n -2})\, (1 - Dq^{2n})},
\ee
\bea
b_n \to \frac{Dq^{2n-3}\, (1 - B^2q^{2n-2})\, (1 - D^2q^{2n
-2}/B^2)}
{B^2\,(1-Dq^{2n-1})\,(1 - Dq^{2n-2})^2\, (1 - D q^{2n-3})}.
\nonumber
\eea
We replace $q$ by $q^{1/2}$ then identify the parameters $A, B, C,
D$ as
\be
B = q^{1 + \a/2} = -C, \quad D = q^{2 + (\a+\b)/2}, \quad A \to
\infty.
\ee
It is not difficult to see that if $Z_n(x)$ satisfies 
\bea
Z_{n+1}(x) = (x + a'_n)\, Z_n(x) \, - b'_n\, Z_{n-1}(x),
\nonumber
\eea
with
\be
a'_n = -\frac{q^{(n-1)/2} (1 + q^{n+(\a + \b +3)/2})\, (1 - q^{(\b
- \a)/2})}
{(1 - q^{n +1+(a+\b)/2})\, (1 - q^{n + 2 +(a+\b)/2})},
\ee
\bea
b'_n = \frac{q^{n+(\b - \a-3)/2}\, (1 - q^{n+\a+1})\, (1 - q^{n
+\b+1})}
{(1-q^{n+(\a+\b+1)/2})\,(1 - q^{n+1+(\a+\b)/2})^2\, (1 - 
q^{n+(\a+\b+3)/2})}.
\nonumber
\eea
It then follows that 
\be
Y_n(x) := q^{(2\a+5)n/4}\,Z_n(xq^{-(2\a+5)/4})
\ee
satisfies (3.12) with $\mu$ replaced by $x$. This relationship
between solutions of (8.8) and (3.12) enables us to take advantage
of the detailed study of solutions of (8.5) contained in
\cite{Gu:Is}. For example one can obtain the minimal solution to
(3.12) by inserting the values of $A, B, C, D$ of (8.8) into
$X_n^{(5)}$
of \cite{Gu:Is}, which remains a minimal solution. This gives an
alternate 
derivation of the form of $Y_k^{(\a, \b)}(x)$ of (5.11).

\bigskip
University of South Florida, Tampa, Florida, 33620, USA.

Carleton University, Ottawa, Ontario, Canada K1S 5B6.

University of Toronto, Toronto, Ontario, Canada M5S 1A1
\end{document}